\documentclass{amsart}
\usepackage{amssymb,amscd, amsmath, amsfonts, enumerate}

\usepackage{graphicx}
\usepackage{color}

\DeclareMathAlphabet{\mathssbx}{OT1}{cmss}{bx}{n}

\begin{document}

\setcounter{tocdepth}{1}

\def\mfd{manifold}
\def\mfds{manifolds}
\def\iff{if, and only if,\ }
\def\g0{G^{(0)}}
\def\h0{H^{(0)}}
\def\maJ{\mathcal{J}}
\def\maK{\mathcal{K}}
\def\codim{\mathrm{codim}}

\renewcommand{\theenumi}{\alph{enumi}}
\renewcommand{\labelenumi}{\rm {({\theenumi})}}
\renewcommand{\labelenumii}{(\roman{enumii})}
%%%%%%%%%%%%%%%%%%%%%%%%%%%%%%%%%%%%%%%%%%%%%%%%%%
%
% OPERATORS:
%
\newcommand{\wt}[1]{\widetilde{#1}}
\newcommand{\Prod}{\prod}
\newcommand{\Cal}{\mathcal}
\newcommand\Ran{\operatorname{Im}}
\newcommand\abs{\operatorname{abs}}
\newcommand\rel{\operatorname{rel}}
\newcommand\inv{\operatorname{inv}}
\newcommand\topo{\operatorname{top}}
\newcommand\opp{\operatorname{op}}
\newcommand\mfk{\mathfrak}
\newcommand\coker{\operatorname{coker}}
\newcommand\hotimes{\hat \otimes}
\newcommand\ind{\operatorname{ind}}
\newcommand\End{\operatorname{End}}
\newcommand\per{\operatorname{per}}
\newcommand\pa{\partial}
\newcommand\sign{\operatorname{sign}}
\newcommand\supp{\operatorname{supp}}
\newcommand\cy{\mathcal{C}^\infty}
\newcommand\CI{\mathcal{C}^\infty}
\newcommand\CO{\mathcal{C}_0}
\newcommand\lra{\longrightarrow}
\newcommand\vlra{-\!\!\!-\!\!\!-\!\!\!\!\longrightarrow}
\newcommand\bS{{}^b\kern-1pt S}
\newcommand\bT{{}^b\kern-1pt T}
\newcommand\Hom{\operatorname{Hom}}

\newcommand\alg[1]{\mathfrak{A}(#1)}
\newcommand\qalg[1]{\mathfrak{B}(#1)}
\newcommand\ralg[1]{\mathfrak{A}_r(#1)}
\newcommand\rqalg[2]{\mathfrak{B}_r(#1)}
\newcommand\ideal[1]{C^*(#1)}
\newcommand\rideal[1]{C^*_r(#1)}
\newcommand\qideal[2]{\mathfrak{R}_{#1}(#2)}
\newcommand\In{\operatorname{In}}

\newcommand\TR{\operatorname{T}}
\newcommand\ha{\frac12}
\newcommand\cal{\mathcal}
\newcommand\END{\operatorname{END}}
\newcommand\ENDG{\END_{\GR}(E)}
\newcommand\CC{\mathbb C}
\newcommand\NN{\mathbb N}
\newcommand\RR{\mathbb R}
\newcommand\R{\mathbb R}
\newcommand\ZZ{\mathbb Z}
\newcommand\ci{${\mathcal C}^{\infty}$}
\newcommand\CIc{{\mathcal C}^{\infty}_{\text{c}}}
\newcommand\hden{{\Omega^{\lambda}_d}}
\newcommand\VD{{\mathcal D}}
\newcommand\mhden{{\Omega^{-1/2}_d}}
\newcommand\ehden{r^*(E)\otimes {\Omega^{\lambda}_d}}

%%%%%%%%%%%%%%%%%%%%%%%%%%%%%%%%%%%%%%%%%%%%%%%%
%%%%% Macros for this manuscript

\newcommand{\calg}{$C^*$-algebra}
\newcommand{\Cat}{\mathcal C}
\newcommand{\Gr}[1]{{\mathcal G}^{(#1)}}
\newcommand{\GR}{\mathcal G}
\newcommand{\LGR}{\mathcal L}
\newcommand{\BB}{\mathbb{B}}
\newcommand{\GG}{\mathcal G}
\newcommand{\OA}{\mathcal O}
\newcommand{\tPS}[1]{\Psi^{#1,0}(\GR)}
%\newcommand{\ttPS}[1]{\Psi_{\loc}^{#1}(\GR;E)}
%\newcommand{\tttPS}[1]{\Psi_{\loc}^{#1}(\GR)}
%% `loc'?
\newcommand{\AL}{{\mathcal A}(\GR)}
\newcommand{\FAM}{P=(P_x,x \in \Gr0)}
\newcommand\symb[2]{{\mathcal S}^{#1}(#2)}
\newcommand{\loc}{\operatorname{loc}}
\newcommand{\cl}{\operatorname{cl}}
\newcommand{\A}{s}
\newcommand{\prop}{\operatorname{prop}}
\newcommand{\comp}{\operatorname{comp}}
\newcommand{\adb}{\operatorname{adb}}
\newcommand{\dist}{\operatorname{dist}}

\newcommand\tgt[1]{{}^{T}\kern-1pt #1}
\newcommand\adi[1]{{}^{ad}\kern-1pt #1}
\newcommand{\alp}{r }
\newcommand{\bet}{d }
\newcommand{\gm}{\Gamma }
\newcommand{\lon}{\longrightarrow }
\newcommand{\be}{\begin{eqnarray*}}
\newcommand{\ee}{\end{eqnarray*}}
\newcommand{\GGR}{{\GR}}
\newcommand{\cald}{{\cal D}}
\newcommand{\calx}{{\cal X}}
\def\cinfo{{\mathcal C}^{\infty,0}}
\def\ccinfo{{\mathcal C}_c^{\infty,0}}
\def\ccinf0{{\mathcal C}_c^{\infty,0}}
\newcommand{\II}{\ccinfo(S^*(\GR), \End(E) \otimes {\mathcal P}_m)}

\newcommand{\III}{\ccinfo(S^*(\GR), \End(E) \otimes {\mathcal P}_m)}

\newcommand{\IIY}{\CIc(S^*(\GR\vert_Y),
\End(E\vert_Y) \otimes {\mathcal P}_m)}
\newcommand{\mI}{\mathfrak I}

\newcommand{\cfg}{Lie groupoid}

\def\nin{\noindent}
\def\eg{e.g.\ }
\def\pt#1#2{{\partial #1\over \partial #2}}

\newcommand{\frakg}{{\mathfrak g}}

\let\Tilde=\widetilde
\let\Bar=\overline
\let\Vec=\overrightarrow
\let\ceV=\overleftarrow
\def\vlra{\hbox{$\,-\!\!\!-\!\!\!-\!\!\!-\!\!\!-\!\!\!
-\!\!\!-\!\!\!-\!\!\!-\!\!\!-\!\!\!\longrightarrow\,$}}

\def\vleq{\hbox{$\,=\!\!\!=\!\!\!=\!\!\!=\!\!\!=\!\!\!
=\!\!\!=\!\!\!=\!\!\!=\!\!\!=\!\!\!=\!\!\!=\!\!\!=\!\!\!=\,$}}

\def\lrah{\hbox{$\,-\!\!\!-\!\!\!
-\!\!\!-\!\!\!-\!\!\!-\!\!\!-\!\!\!\longrightarrow\,$}}

\def\surj{-\!\!\!-\!\!\!-\!\!\!\gg}

\def\inj{>\!\!\!-\!\!\!-\!\!\!-\!\!\!>}

\newcommand{\ad}[1]{{\mathcal A}^{\fd}_{#1}}
\newcommand{\as}{\wt{a}}
\newcommand{\aun}{\mathcal{A}_{\infty}}
\newcommand{\aund}{\ad{\infty}}
\newcommand{\bo}{{}^{b}\Omega}
\newcommand{\bu}{{\mathcal B}_{U}}
\newcommand{\bs}{\wt{b}}
\newcommand{\ca}{{\mathcal A}}
\newcommand{\cae}{\ca_{1}}
\newcommand{\can}{\ca_{0}}
\newcommand{\caun}{\ca_{\infty}}
\newcommand{\cav}{\ca_{\Omega}}
\newcommand{\cavs}{\wt{\ca}_{\Omega}}
\newcommand{\maC}{\mathcal C}
\newcommand{\cb}{{\mathcal B}}
\newcommand{\cH}{{\mathcal{H}}}
\newcommand{\cj}{{\mathcal J}}
\newcommand{\cjt}[1]{\cj_{#1}(\ft)}
\newcommand{\cjtk}{\cj_{k}(\ft)}
\newcommand{\cjtu}{\cjt{\infty}}
\newcommand{\ck}{{\mathcal K}}
\newcommand{\ckx}{{\mathcal K}_{x}}
\newcommand{\cnb}{\mathcal{C}_{b}}
\newcommand{\cun}{\mathcal{C}^{\infty}}
\newcommand{\cunb}{\mathcal{C}^{\infty}_{b}}
\newcommand{\cunc}{\cun_{c}}
\newcommand{\cuncd}{\dot{\mathcal{C}}^{\infty}_{c}}
\newcommand{\cv}{\mathcal{V}}
\newcommand{\cz}{\CC}
\newcommand{\dbar}{d\hspace{-2.8pt}\rule[5.5pt]{3pt}{0.24pt}}
\newcommand{\dd}{\Delta_{\fd}}
\newcommand{\de}[1]{{\mathcal D}(#1)}
\newcommand{\deh}{{\mathcal D}^{\frac{1}{2}}}
\newcommand{\dez}{\delta_{12}}
\newcommand{\dirint}{\int^{\oplus}}
\newcommand{\djx}{D_{j;x}}
\newcommand{\djxs}{\wt{D}_{j;x}}
\newcommand{\dmy}{D_{m;y}}
\newcommand{\djys}{\wt{D}_{j;y}}
\newcommand{\dmjx}{\delta_{M_{j,x}}}
\newcommand{\dt}{\delta_{T}}
\newcommand{\dxs}{\frac{dx'}{x'}}
\newcommand{\dxss}{dx''}
\newcommand{\elm}{E_{\ell,m}}
\newcommand{\fd}{\mathfrak{D}}
\newcommand{\fm}{\mathfrak{M}}
\newcommand{\fs}{\mathfrak{S}}
\newcommand{\ft}{\mathfrak{T}}
\newcommand{\gpt}{\mathcal{G}_{\Psi}(T)}
\newcommand{\gx}{{\GR_{x}}}
\newcommand{\hol}[1]{{\cal O}(#1)}
\newcommand{\hft}[1]{{\mathcal H}^{#1}(\ft)}
\newcommand{\huft}{\hft{\infty}}
\newcommand{\idh}{{\rm id}_{H}}
\newcommand{\lch}{\mathcal{L}(\cH)}
\newcommand{\lck}{\mathcal{L}(\ck)}
\newcommand{\lckx}{\mathcal{L}(\ckx)}
\newcommand{\ldch}{\mathcal{L}_{D}(\cH)}
\newcommand{\ldgch}{\mathcal{L}^{\GR}_{D}(\cH)}
\newcommand{\lh}{\mathcal{L}(H)}
\newcommand{\lhx}{\mathcal{L}(H_{x})}
\newcommand{\llzrny}{\mathcal{L}(\lzrny)}
\newcommand{\lukdh}{{\mathcal L}^{\infty}_{K,d}(H)}
\newcommand{\lukh}{{\mathcal L}^{\infty}_{K}(H)}
\newcommand{\lzrny}{L^{2}(\rny)}
\newcommand{\lzgx}{L^{2}(\gx;\rsed)}
\newcommand{\mjx}{M_{j;x}}
\newcommand{\mmy}{M_{m;y}}
\newcommand{\norm}[1]{\|#1\|}
\newcommand{\nz}{\NN}
\newcommand{\olt}{\omega_{T}^{\ell}}
\newcommand{\olta}[1]{\omega_{#1}^{\ell}}
\newcommand{\olrt}{\omega_{T_{1},T_{2}}^{\ell,r}}
\newcommand{\ort}{\omega_{T}^{r}}
\newcommand{\orta}[1]{\omega_{#1}^{r}}
\newcommand{\ot}{\omega_{T}}
\newcommand{\otimesh}{\wh{\otimes}}
\newcommand{\pege}{\Psi^{1,0}(\GR;\deh)}
\newcommand{\pft}[1]{\Psi^{#1}(\ft)}
\newcommand{\pmge}{\Psi^{m,0}(\GR;\deh)}
\newcommand{\pmunge}{\Psi^{-\infty,0}(\GR)}
\newcommand{\pnge}{\Psi^{0,0}(\GR)}
\newcommand{\pmungee}{\Psi^{-\infty,0}(\GR)}
\newcommand{\pngee}{\Psi^{0,0}(\GR)}
\newcommand{\puft}{\pft{\infty}}
\newcommand{\re}{R^{(1)}}
\newcommand{\rklx}{\RR^{k,\ell}_{x}}
\newcommand{\rkx}{\RR^{k}_{x}}
\newcommand{\rne}{\RR^{n}_{\eta}}
\newcommand{\rny}{\RR^{n}_{y}}
\newcommand{\rpq}{\overline{\rz}_{+}}
\newcommand{\rsed}{r^{*}\deh}
\newcommand{\rz}{\RR}
\newcommand{\rj}{R^{(j)}}
\newcommand{\sk}{\Sigma_{K}}
\newcommand{\spk}[1]{\langle #1 \rangle}
\newcommand{\sva}{\sigma_{\Omega}(a;\varphi_{0},\psi_{0})}
\newcommand{\te}{T^{(1)}}
\newcommand{\tz}{T^{(2)}}
\newcommand{\tj}{T^{(j)}}
\newcommand{\wh}[1]{\widehat{#1}}

\newcommand\ie{{\em i.e.,\ }}
%
% THEOREM TYPE ENVIRONMENTS:
%
\newtheorem{theorem}{Theorem}
\newtheorem{proposition}{Proposition}
\newtheorem{corollary}{Corollary}
\newtheorem{lemma}{Lemma}
\newtheorem{definition}{Definition}
\newtheorem{notation}{Notations}
\theoremstyle{remark}
\newtheorem{remark}[theorem]{Remark}
\newtheorem{example}[theorem]{Example}
\newtheorem{examples}[theorem]{Examples}

%

%_
%%%%%%%%%%%%%%%%%%%%%%%%%%%%%%%%%
%%                              %
%%  MACROS FOR THIS MANUSCRIPT  %
%%                              %
%%%%%%%%%%%%%%%%%%%%%%%%%%%%%%%%%
%
% NEW, FOR THIS MANUSCRIPT
%
% \newcommand{\cb}{{\mathcal B}}
% \newcommand{\ck}{{\mathcal K}}
\def\n#1#2{\| #1 \|_{{#2}}}
\def\g0{\GR^{(0)}}
\def\S{\mathscr{S}}
\def\cred#1{C^*_{\mathrm{r}}(#1)}
\def\rep#1{{\R_+^*}\!^{#1}}  %%R Etoile Plus
\def\r+#1{\R_+^{#1}}
\def\cinfo{{\mathcal C}^{\infty,0}}
\def\ccinfo{{\mathcal C}_c^{\infty,0}}
\def\Cc{{\mathcal C}_c}
%%%%%%%%%%%%%%%%%%%%%%%%%%%%%%%%%
%%                              %
%%  MACROS B.M.                 %
%%                              %
%%%%%%%%%%%%%%%%%%%%%%%%%%%%%%%%%
\let \mx \mbox
\let \hx \hbox
\let \vx \vbox

\def\C{\mathbb{C}}%c'est eclaire
\def\E{\mathbb{E}}
\def\A{\mathbb{A}}
\def\B{\mathbb{B}}
\def\D{\mathbb{D}}
\def\F{\mathbb{F}}
\def\H{\mathbb{H}}
\def\N{\mathbb{N}}
\def\P{\mathbb{P}}
\def\Q{\mathbb{Q}}
\def\R{\mathbb{R}}
\def\Z{\mathbb{Z}}
\def\G{\mathbb{G}}
\def\K{\mathbb{K}}
\def\J{\mathbb{J}}
\def\T{\mathbb{T}}

\def\diagr#1{\def\normalbaselines{\baselineskip=0pt
\lineskip=10pt\lineskiplimit=1pt} \matrix{#1}}
\def\hfl#1#2{\smash{\mathop{\hx to 12mm{\rightarrowfill}}
\limits^{\scriptstyle#1}_{\scriptstyle#2}}}
\def\vfl#1#2{\llap{$\scriptstyle #1$}\left\downarrow
\vx to 6mm{}\right.\rlap{$\scriptstyle #2$}}
\def\antihfl#1#2{\smash{\mathop{\hx to 12mm{\leftarrowfill}}
\limits^{\scriptstyle#1}_{\scriptstyle#2}}}
\def\antivfl#1#2{\llap{$\scriptstyle #1$}\left\uparrow
\vx to 6mm{}\right.\rlap{$\scriptstyle #2$}}

\def\build#1_#2^#3{\mathrel{\mathop{\kern 0pt#1}\limits_{#2}^{#3}}}

\def\limind{\mathop{\oalign{lim\cr\hidewidth$\longrightarrow$\hidewidth\cr}}}
\def\limproj{\mathop{\oalign{lim\cr \hidewidth$\longleftarrow$\hidewidth\cr}}}

\def\norme#1{\left\| #1 \right\|}
\def\module#1{\left| #1 \right|}
\def\va#1{\left| #1 \right|}
\def\scal#1{\left\langle #1 \right\rangle}
\def\inv#1{{#1}^{-1}}
\def\invf{f^{-1}}

\def\g0{G^{(0)}}
\def\psd{pseudodiff\'erentiel}
\def\dg{\partial G}
\def\dm{\partial M}
\def\dx{\partial X}

\def\intm{\overset{\:\bullet}{M}}
\def\intx{\overset{\:\bullet}{X}}
\def\intf#1{\overset{\:\bullet}{#1}}
\def\cc#1{C^*(#1)}
\def\kc#1{K_*(C^*(#1))}
\def\rep#1{{\R_+^*}^{#1}}  
\def\r+#1{\R_+^{#1}}

\def\tg{groupo\"\i de tangent}
\def\gr{groupoid}
\def\d{{\rm d}}
\def\e{{\varepsilon}}
\def\cstar{$C^*$-algebra}
\def\rbar{\overline\R}
\def\l{\lambda}
\def\kth{$K$-theory}
\def\cinf{$C^\infty$}
\def\M{M}
\def\G{{\bf G}}
\def\E{\mathcal{E}}
\def\In{\hbox{In}}
\def\pd{pseudodifferential}

\def\k{\mathssbx{k}}

\def\lb{\hbox{lb}}
\def\rb{\hbox{rb}}

\def\bsp{$b$-stretched product}
\def\bcalc{$b$-calculus}

\def\sm{submanifold}
\def\mc{manifold with corners}
\def\mecs{manifolds with embedded faces}
\def\mec{manifold with embedded faces}
\def\mcs{manifolds with corners}

\def\D{\hbox{D}}
\def\S{\mathcal S}

\def\A{\mathcal A}
\def\F{\mathcal F}
\def\HF{\mathcal{HF}}
\def\sym{\mathfrak S}
\def\GG{\mathcal G}
\def\sh{\mathrm{sh}}
\def\ch{\mathrm{ch}}

\def\fr{\frac}
\def\ub{\underbar}
\def\O{\Cal O}
\def\F{\Cal F}
\def\differ{\text{differentiable} }
\def\tPSeudo{\text{pseudodifferential} }
\def\supp{\text{supp} }
\def\inn{{\mathcal R}}
\def\prop{\text{prop}}
\def\frag{\frak{G}}
\def\simd{\tilde{d}}
\def\simf{\tilde{\F}}
\def\simo{\tilde{\O}}
\def\simr{\tilde{r}}
\def\simp{\tilde{p}}
\def\simmu{\tilde{\mu}}
\def\O12{\Omega^\frac{1}{2}}

\def\maG{\mathcal{G}}
\def\maH{\mathcal{H}}
\def\maK{\mathcal{K}}
\def\maI{\mathcal{I}}
\def\maL{\mathcal{L}}
\def\maV{\mathcal{V}}
\newcommand{\Diff}{\operatorname{Diff}}
\newcommand{\Diffb}{\Diff_b(M)}
\newcommand{\DiffVM}{\operatorname{\Diff_{\maV}(M)}}
\newcommand{\PsiVM}{\Psi_\maV^\infty(M)}
\newcommand{\Psib}{\Psi_b^\infty(M)}

\date\today
\author[B. Monthubert]{Bertrand Monthubert}
       \address{ Universit\'e Paul
       Sabatier (UFR MIG), Laboratoire Emile Picard,
       F-31062 Toulouse CEDEX 4}
       \email{bertrand.monthubert@math.ups-tlse.fr}

\author[V. Nistor]{Victor Nistor} \address{Pennsylvania State
       University, Math. Dept., University Park, PA 16802}
       \email{nistor@math.psu.edu}

\thanks{Monthubert was partially supported by a ACI Jeunes
       Chercheurs. Manuscripts available from {\bf
       http:{\scriptsize//}bertrand.monthubert.net}.
       Nistor was partially supported by the NSF Grant DMS
       0200808. Manuscripts available from {\bf
http:{\scriptsize//}www.math.psu.edu{\scriptsize/}nistor{\scriptsize/}}.}

%\dedicatory\datverp
\begin{abstract}
We define an analytic index and prove a topological index theorem
for a non-compact manifold $M_0$ with poly-cylindrical ends. We
prove that an elliptic operator $P$ on $M_0$ has an invertible
perturbation $P+R$ by a lower order operator \iff its analytic
index vanishes. As an application, we determine the $K$-theory
groups of groupoid $C^*$--algebras of manifolds with corners.
\end{abstract}

\title[An index theorem for manifolds with corners]{A topological
index theorem for manifolds with corners}

\maketitle \tableofcontents

\section*{Introduction}

Let $M_0$ be a smooth, {\em compact} manifold and $D$ be an
elliptic differential operator of order $m$ acting between smooth
sections of vector bundles on $M$. Then $D$ is continuous and
Fredholm as a map between suitable Sobolev spaces (\ie $H^s \to
H^{s-m}$). In particular, the kernel and cokernel of $D$ are
finite dimensional. This allows one to define $\ind(D)$, the {\em
Fredholm index} of $D$, by
\begin{equation*}
    \ind(D) := \dim \ker(D) - \dim \coker(D) =
    \dim \ker(D) - \dim (H^{s-m}/DH^s).
\end{equation*}
The knowledge of the Fredholm index is relevant because it gives
an obstruction to $D$ being invertible. More precisely, we have
the following result whose proof is an easy exercise in Functional
Analysis.

\begin{theorem}\label{thm.perturb}\
There exists a pseudodifferential operator $R$ of order $< m =
\operatorname{ord}(D)$ on the compact manifold $M_0$ such that $D
+ R : H^s \to H^{s-m}$ is an isomorphism \iff $\ind(D) = 0$.
Moreover, if one such $R$ exists, then it can be chosen to be of
order $-\infty$, i.e., regularizing.
\end{theorem}

If $M_0$ is non-compact, the case of main interest in our paper,
then an elliptic differential operator $D$ on $M_0$ needs not be
Fredholm in general (here, by ``elliptic'' we mean that the
principal symbol is invertible). The Fredholm index of $D$ is
therefore not defined and Theorem \ref{thm.perturb} is meaningless
as it stands. We therefore introduce an extension of the Fredholm
index on manifolds with poly-cylindrical ends. A {\em manifold
with poly-cylindrical ends} is locally diffeomorphic at infinity
with a product of manifolds with cylindrical ends (see Definition
\ref{def.pce}).

{\em From now on, $M_0$ will denote a manifold with
poly-cylindrical ends, usually non-compact, and $M$ will denote
its given compactification to a manifold with corners.} Two of the
main results of this paper are an extension of the Atiyah--Singer
Index Theorem on the equality of the topological and analytic
index and an extension of Theorem \ref{thm.perturb} to manifolds
with poly-cylindrical ends (with the analytic index replacing the
Fredholm index). Our work is motivated, in part, by recent work of
Leichtnam--Piazza \cite{LP1} and Nistor--Troitsky \cite{NT1},
who have showed that suitable generalizations of Theorem
\ref{thm.perturb} are useful in geometric and topological
applications. Extensions of Theorem \ref{thm.perturb} seem to be
important for the study of boundary value problems on polyhedral
domains \cite{N-prden} using the method of layer potentials. In
the process, we establish several extensions to non-compact
manifolds with poly-cylindrical ends of the classical properties
on the topological index of elliptic operators on compact
manifolds \cite{AS1, carvalho}.

Let us explain in a little more detail our main results. The
tangent bundle $TM_0$ extends to a bundle $A_M$ on $M$ whose
sections are the vector fields tangent to the faces of $M$ (this
is the ``compressed tangent bundle'' in Melrose's terminology).
Let $D$ be a differential (or pseudodifferential operator) on
$M_0$ compatible with the structure at infinity on $M_0$ (\ie $P
\in \Psi^{\infty}_b(M)$, where $\Psi^{\infty}_b(M)$ denotes the
algebra of pseudodifferential operators on $M_0$ compatible with
the structure at infinity). Then the principal symbol of $D$
extends to a symbol defined on $A_M^*$.

Assume that $D$ is elliptic, in the sense that its principal
symbol is invertible on $A_M^*$ outside the zero section. Then $D$
is invertible modulo regularizing operators. Let $C^*(M)$ be the
norm closure of the algebra of regularizing operators on $M_0$
compatible with the structure at infinity of $M_0$. The
obstruction of $D$ to be invertible defines a map
\begin{equation}
    \ind_a = \ind_a^M : K^0(A_M^*) \to K_0(C^*(M)).
\end{equation}
This extension of the topological index is the {\em analytic
index} map mentioned above. In particular, the principal symbol of
$D$ defines an element $\ind_a(D) \in K_0(C^*(M))$. Theorem
\ref{thm.perturb} then remains true for manifolds with
poly-cylindrical ends if we use the analyitic in place of the
Fredholm index.

To explain our generalization of the topological index theorem,
consider an embedding $\iota : M \to X$ of manifolds with corners
and let $\iota_{!}$ be the push-forward map in $K$-theory. Then
Theorem \ref{commutativity} states the commutativity of the
diagram
\begin{equation}\label{diag.I}
\begin{CD}
    K_0(C^*(M)) @>{\iota_*}>> K_0(C^*(X))\\
    @A{\ind_a^M}AA  @AA{\ind_a^X}A\\
    K^0(A_M^*) @>{\iota_{!}}>>  K^0(A_X^*),\\
\end{CD}
\end{equation}

In order to interpret the above diagram as an index theorem in the
usual way, we need to choose $X$ such that $\iota_* : K_0(C^*(M))
\to K_0(C^*(X))$ and $\ind_a^X : K^0(A_X^*) \to K_0(C^*(X))$ be
isomorphisms. This would the provide us with both an
identification of the groups $K_0(C^*(M))$ and of the map
$\ind_a^M$.

To obtain a manifold $X$ with the above mentioned properties, we
shall proceed as follows. Let us assume $M$ is compact with
embedded faces (recall \cite{MelroseScattering} that this
means that each face of $M$ of maximal dimension has a defining
function). To $M$ we will associate a non-canonical \mec\ $X_M$
and an embedding $i : M \to X_M$ such that
\begin{enumerate}[(i)]
\item\ each open face of $X_M$ is diffeomorphic to a Euclidean
space,
\item\ each face of $M$ is the transverse intersection of $M$ and
of a face of $X_M$,
\item\ $F \to F \cap M$ establishes a bijection between the open
faces of $M$ and those of $X_M$.
\end{enumerate}
We shall say that an embedding $X_M$ with the above properties is a
classifying space for $M$ and prove that $i_K : K_0(C^*(M)) \to
K_0(C^*(X_M))$ is an isomorphism, where $j_K$ is a canonical morphism
associated to any embedding of manifolds with corners $j$ (see Lemma
\ref{lemma.prop.3}).

Let us now summarize the contents of the paper. In Section
\ref{basic}, we review the definitions of manifolds with corners, of
Lie groupoids, and of Lie algebroids. We also review and extend a
result from \cite{nistorINT} on the integration of Lie algebroids. In
Section \ref{a.index} we recall the definition of the analytic index
using the tangent groupoid and then compare this definition with other
possible definitions of an analytic index. In the process, we
establish several technical results on tangent groupoids. As an
application, in Theorem \ref{thm.perturb.M}, we provide the
generalization of Theorem \ref{thm.perturb} mentioned above.  Section
\ref{sec.prop} contains the main properties of the analytic index. In
this section, we also introduce the morphism $j_K$ associated to an
embedding of manifolds with corners $j$ and we provide conditions for
$j_K$ and $\ind_a^M$ to be isomorphisms.  The compatibility of the
analytic index and of the shriek maps is established in Section
\ref{commutativity-diagram}. This is then used in the following
section to establish the equality of the analytic and topological
index. In the last section we show that it is indeed possible to find
a space $X_M$ and an embedding $i : M \to X_M$ with the properties
(i)--(iii) above (\ie a classifying space for $M$ exists).

\thanks{We thank Catarina Carvalho, Sergiu Moroianu, and Georges
Skandalis for useful discussions. The first named author also
wishes to acknowledge the constant support and encouragement of
his former Ph.D. advisor, Georges Skandalis.}

\section{Basic definitions}\label{basic}

In this section we recall several basic definitions and
constructions, including:\ manifolds with corners, manifolds with
poly-cylindrical ends, groupoids, adapted pseudodifferential
operators, and Lie algebroids. Good references  to the special
issues we deal with in this section are \cite{bm1, bm-jfa, LMN1},
and \cite{NWX}. General references to the subject are
\cite{ConnesNCG, Mackenzie, MelroseScattering}, and
\cite{Renault1}.

\subsection{Manifolds with corners}
A manifold with corners $M$ is a manifold modeled on $(\R_+)^n$
(which is denoted $\R_+^n$, and must not be understood as a
half-space of $\R^n$ but as a quadrant). This means that any point
$x \in M$ has a neighborhood of the form $\R^{n-k}\times \R_+^k$
(with $x$ mapping to $0 = (0, 0, \ldots, 0)$). We call $k$ the
\textit{depth} of $x$.

The set of points of depth $k$ is a union of connected components,
each of which is a smooth submanifold in its own, called an
\textit{open face of codimension $k$}. The closure of an open face
(of codimension $k$) is just called a \textit{face} of $M$ (of
codimension $k$). A closed face is not necessarily a
manifold with corners (think of the ``tear-drop domain'' in the
plane). A codimension one face of $M$ is called a
\textit{hyperface} of $M$.

We shall sometimes require that each hyperface of $M$ be an
embedded submanifold of $M$. If this is the case, we shall call
$M$ a \textit{\mec}, while in Melrose's terminology it is just
called a \textit{\mc}. A manifold $M$ is a \mec\ if each hyperface
has a \textit{defining function}. Recall that a function $\rho : M
\to \R_+$ is a defining function of the hyperface $H \subset M$ if
$\rho$ is smooth, if $\rho(x)=0$ precisely when $x \in H$, and
$\d\rho$ does not vanish on $H$. Such a defining function provides
us, in particular, with a trivialization of $NH = TM/TH$, the
normal bundle of $H$.

If $X$ and $Y$ are two \mcs, a map $f:X \to Y$ is called a {\em
closed embedding} of \mcs\ if
\begin{enumerate}
\item\ it is differentiable, injective, with closed range;
\item\ $df$ is injective;
\item\ for each open face $F$ of $Y$, $f(X)$ is transverse to $F$
(recall that this means that $df(T_xX) + T_yF = T_yY$, if $y =
f(x) \in F$); and
\item\ each hyperface of $X$ is a connected component of the
inverse image of a hyperface of $Y$.
\end{enumerate}
In particular, $x$ and $f(x)$ will have the same depth.

A {\em submersion} $f : M \to N$ (between two manifolds with
corners $M$ and $N$) is a differentiable map $f$ such that
\begin{enumerate}
\item\ $df(v)$ is an inward pointing tangent vector of $N$ if, and
only if, $v$ is an inward pointing vector of $M$; and
\item\ $df$ is surjective at all points.
\end{enumerate}

If $\iota:X \to Y$ is a closed embedding, then there exists a
\textit{tubular neighborhood} $U \subset Y$, which means that
there exists a vector bundle $E$ over $X$, and an open
neighborhood $Z$ of the zero section isomorphic to $U$, such that
the following diagram is commutative:
\begin{equation*}
\begin{CD}
    Z @>\simeq>> U\\ @AAA @VVV\\ X @>i>> Y\\
\end{CD}
\end{equation*}

The existence of such a tubular neighborhood is proved in
\cite{douady-sem-cartan, AIN}. For the benefit of the reader, we
now provide a sketch of the proof formulated in our framework. Let
us denote by $\maV_b(X) \subset \Gamma(TX)$ the subspace of vector
fields {\em tangent} to all faces of $X$. Then $\maV_b(X)$ is a
Lie algebra with respect to the Lie bracket and a
$\CI(X)$--module. Moreover, there exists a vector bundle $A_X \to
X$, uniquely determined up to isomorphism, such that
\begin{equation*}
    \maV_b(X) \simeq \Gamma(A_X) \quad \text{as }\;
    \CI(X) \text{--modules}.
\end{equation*}
If $Y \subset X$ is a closed, embedded submanifold with
corners, then $\maV_b(Y)$ consists of the restrictions to $Y$ of
the vector fields $V \in \maV_b(M)$ with the property that $V$ is
tangent to $Y$. Fix an arbitrary open face $F$ of $Y$. The
assumption that every hyperface of $Y$ be contained in a hyperface
of $X$ (even a connected component of the intersection of $Y$ with
a hyperface of $X$) implies that there exists an open face $F'$ of
$X$, of the same codimension as $F$, such that $F \subset F'$.
Next, the definition of an embedded submanifold with corners
implies that $Y$ is transverse to $F'$ and hence the natural map
$T_xF'/T_xF \to T_xX/T_xY$ is an isomorphism for any $x \in F$. It
is then possible\footnote{This was done in the general framework
of Lie manifolds in \cite{aln1} and was used in \cite{AIN}.} to
construct a connection whose associated exponential map
$exp:N_X^Y\to Y$ is such that its restrictions to each face is a
map $exp_{|F}:(N_X^Y)_F\to F'$ and is a local diffeomorphism
around the zero section. Using this exponential map, one can
define $E$, $Z$, and $U$ with the required properties.

\subsection{Manifolds with poly-cylindrical ends} We shall use the
notation and terminology introduced in the previous subsection.
Let us fix a metric on the bundle $A_M \to M$. This metric
defines, in particular, a Riemannian metric on $M_0$, the interior
of $M$. The following definition is from \cite{MelroseScattering}

\begin{definition}\label{def.pce}\
A manifold with poly-cylindrical ends is a smooth Riemannian
manifold $M_0$ that is diffeomorphic with the interior of a
compact manifold with corners $M$ such that the metric on $M_0$ is
the restriction of a metric on $A_M$.
\end{definition}

The interior $M_0$ of any compact manifold with corners $M$ is
therefore a manifold with poly-cylindrical ends, for any choice of
a metric on $A_M \to M$. The compact manifold with corners $M$
will be called the {\em compactification of $M_0$}.

Let us denote the hyperfaces of $M$ by $H_i$, $i = 1, \ldots, N$,
and fix, for any hyperface $H_i$, a defining function $\rho_i$, if
that hyperface has such a function (\ie\ if $H_i$ is an embedded
hyperface). Let $h$ be a smooth metric on a manifold with embedded
faces $M$. Then a typical example of a manifold with
poly-cylindrical ends is provided by $M_0$, the interior of $M$,
with the metric
\begin{equation}\label{eq.ex.pce}
    g = h + \sum_{i=1}^N (\rho_i^{-1}d\rho_i)^2.
\end{equation}
For instance, if we take $M = [0, 1]^n$ with defining functions
$x_j$ and $1 - x_j$ and let $h = 0$, then the resulting manifold
with poly-cylindrical ends is isometrically diffeomorphic with
$\R^n$ with the standard (flat) Euclidean metric.

We shall denote by $\Diffb$ the algebra of differential operators
generated by $\maV_b(M)$ together with multiplications by
functions in $\CI(M)$. A simple but useful result \cite{aln1}
states that all geometric operators (Laplace, Dirac, signature,
... ) associated to the Riemannian manifold with poly-cylindrical
ends $M_0$ are in $\Diffb$. We shall therefore restrict our study
of differential operators on $M_0$ to differential operators in
$\Diffb$.

If each face of $M$ has a defining function (\ie {\em $M$ has
embedded faces}), then in \cite{MelrosePiazza, MelroseScattering}
it was constructed an algebra $\Psib$ of pseudodifferential
operators on $M$. One of the main properties of $\Psib$ is that a
differential operator $P$ is in $\Psib$ precisely when $P \in
\Diffb$. This construction was generalized in \cite{bm1} to
arbitrary manifolds with corners. For the purpose of this paper,
it is convenient to introduce the algebra $\Psib$ and related
algebras using Lie groupoids. (In fact, our algebras are slightly
smaller than the ones in \cite{MelrosePiazza, MelroseScattering}.)
An operator $P \in \Psib$ will be called {\em compatible with the
structure at infinity} on $M_0$.

\subsection{Differentiable groupoids}
A {\em small category} is a category whose class of
morphisms is a set. The class of objects of a small category is
then a set as well. By definition, a {\em groupoid} is a small
category $\GR$ in which every morphism is invertible. See
\cite{Renault1} for general references on groupoids.

We shall follow the general notations: the set of objects (or {\em
units}) of a groupoid $\GR$ is denoted by $\GR^{(0)}$, and the set
of morphisms (or {\em arrows}) is denoted, by abuse of notation,
by $\GR$ instead of $\GR^{(1)}.$ A \gr\ is endowed with two maps,
the \textit{domain} $d:\GR \to \GR^{(0)}$ and the \textit{range}
$r:\GR \to \GR^{(0)}$. The multiplication $gh$ of
$g,h\in\GR^{(0)}$ is defined if, and only if, $d(g) = r(h)$.  A
groupoid $\GR$ is completely determined by the spaces $\GR^{(0)}$
and $\GR$ and by the structural morphisms: $d$, $r$,
multiplication, inversion, and the inclusion $\GR^{(0)} \to \GR$.

We shall consider {\em Lie groupoids} $(\GR, M)$, that is,
groupoids endowed with a differential structure such that the set
of arrows, $\GR$, and the set of units, $M$, are smooth manifolds
with corners, all structural maps are differentiable, and $d$ is a
submersion of manifolds with corners. In particular, $d^{-1}(x)$
is a smooth manifold (without corners) for any $x \in M$ and
$\GR^{(0)}$ is an embedded submanifold with corners of $\GR$. (The
terminology ``differentiable groupoid'' was used in \cite{NWX}
instead of ``Lie groupoid,'' because the name ``Lie groupoid'' was
used in the past for differentiable groupoids with additional
structures. This has changed, however, and the terminology ``Lie
groupoid'' better reflects the current use.)

A Lie groupoid $\maG$ is called {\em $d$-connected} \iff all the
sets $\maG_x := d^{-1}(x)$ are connected (and hence also path
connected). If we are given a Lie groupoid $\maG$, let us define
$\maG_0$ to consist of all the path components of the units in the
fibers $\maG_x$. Then $\maG_0$ is the an open subset of $\maG$
containing the units and is closed under the groupoid operations.
We shall call $\maG_0$ the $d$-connected component of the units in
$\maG$. It is a Lie groupoid on its own, and, as such, it is
$d$-connected.

\begin{examples}Let us include here some examples of Lie groupoids that will
be needed later on.

\begin{enumerate}[(1)]
\item If $X$ is a smooth manifold, $X \times X$, the {\em pair
groupoid} has units $X$ and is defined by $d(x,y) = y$, $r(x, y) =
y$, and $(x, y) (y, z) = (x, z)$.

\item Let $\pi : X \to M$ be a fibration with smooth fibers, $M$ a
manifold with cornesr. Recall then that the {\em fiberwise product
groupoid} $\GR := X \times_M X$ also has units $X$ and is defined
as $\GR := \{(x_1, x_2) \in X^2, \pi(x_1) = \pi(x_2)\}$ with units
$X$, $d(x_1, x_2) = x_2$, $r(x_1, x_2) = x_1$, and product $(x_1,
x_2) (x_2, x_3) = (x_1, x_3)$. Thus $\GR$ is a subgroupoid of the
pair groupoid $X \times X$. This example will be needed later on
in the proof of Proposition \ref{prop.diag1}.

\item If $\GR_j$, $j=1,2$, are Lie groupoids, then $\GR_1 \times
\GR_2$ is also a Lie groupoid.

\item If $M$ is a manifold with corners, $\GR = M$ is a Lie groupid
with only units.
\end{enumerate}
\end{examples}

\subsection{Lie algebroids}
This subsection may be skipped at a first reading. A \textit{Lie
algebroid} $\pi : A \to M$ over a smooth manifold with corners $M$
is a smooth vector bundle $A$ over $M$ for which there is given a
vector bundle map $\varrho : A \to TM$ satisfying:
\begin{enumerate}
    \item\ $\Gamma(A)$ is endowed with a Lie algebra structure;
    \item\ $[\varrho(X), \varrho(Y)]=\varrho([X,Y])$;
    \item\ $[X, fY]=(\varrho(X)f) Y  + f[X,Y]$;
\end{enumerate}
for any $\CI$--sections $X$ and $Y$ of $A$ and any $\CI$--function
$f$ on $M$.

The simplest example of a Lie algebroid is the tangent bundle $TM
\to M$, with the Lie algebra structure on $\Gamma(TM)$ being given
by the Lie bracket. Similarly, the vector bundle $A_M \to M$
introduced above is also a Lie algebroid, the Lie algebra
structure being again given by the Lie bracket.

If $\GR$ is a Lie groupoid, let us denote by
\begin{equation*}
    T_d\GR := \bigcup_{g\in\GR} T_{g} \GR_{d(g)} \quad \text{and }\,
    A(\GR) := \bigcup_{x \in \GR^{(0)}} T_{x} \GR_{x} =
    T_d\GR\vert_{\GR^{(0)}}
\end{equation*}
the {\em $d$-vertical tangent bundle of $\GR$} and, respectively,
the {\em Lie algebroid of $\GR$}. The sections of $T_d\GR$ are
vector fields tangent to the fibers of $d$, and hence the space of
all these sections, $\Gamma(T_d\GR)$, is closed under the Lie
bracket. Multiplication to the right by $\gamma$ maps $\GR_y$ to
$\GR_x$, where $x$ and $y$ are the domain and the range of $\gamma
\in \GR$. A vector field $X \in \Gamma(T_d\GR)$ will be called
{\rm right invariant} if it is invariant under all these right
multiplications. The space $\Gamma_R(T_d\GR)$ of right invariant
vector fields is also closed under the Lie bracket because a local
diffeomorphism preserves the Lie bracket. Since $A(\GR)$ is the
restriction of the $d$-vertical tangent bundle $T_d\GR$ to the
space of units $M$ and $r$ is a submersion, we have
$\Gamma_R(T_d\GR) \simeq \Gamma(A(\GR))$ and hence the later is a
Lie algebra. It turns out that $A(\GR)$ is, indeed, a Lie
algebroid, where $\varrho:=r_*$, the differential of the map $r$.

In the following, we shall use the following result from
\cite{nistorINT}, formulated in the way it will be used in this
paper. For any set $S$, we shall denote by $S^c$ its complement
(in a larger set that is understood from the context).

\begin{theorem}\label{thm.int}
Let $A \to M$ be a Lie algebroid and anchor map $\varrho : A \to
TM$. Let $N \subset F$ be a closed submanifold (possibly with
corners) of a face $F$ of $M$. Assume that $N$ is invariant, in
the sense that $\varrho(X)$ is tangent to $N$ for any $X \in
\Gamma(A)$. Let $\GR_1$ be a groupoid with units $N$ and Lie
algebroid $A\vert_{N}$. Also, let $\maG_2$ be a groupoid with
units $N^c = M \smallsetminus N$ and Lie algebroid $A\vert_{N^c}$.
Assume that both $\maG_1$ and $\maG_2$ are $d$-connected. Then the
disjoint union $\maG := \maG_1 \cup \maG_2$ has at most one
differentiable structure compatible with the differentiable
structures on $\maG_1$ and $\maG_2$ that makes it a Lie groupoid
with Lie algebroid $A$.
\end{theorem}

Let us remark that by abstract set theory nonsense, it is enough
to require in the above theorem that $A(\GR_1) \simeq A\vert_{N}$
and $A(\GR_2) \simeq A\vert_{N^c},$ but then we have to make sure
that the isomorphism $A(\maG) \simeq A$ is such that it restricts
to the given isomorphisms on $N$ and $N^c$.

The differentiable structure, if there is one, is obtained using
the exponential maps (see \cite{nistorINT} for details). It may
happen, however, that there is no differentiable structure with
the given properties on $\maG$ or that the resulting manifold is
non-Hausdorff. Variants of the above result for stratifications
with finitely many strata can be obtained by induction.

A corollary of the above result is the following. We use the same
notation as in the previous theorem.

\begin{corollary}\label{cor.int}\
Assume $\GR$  has a smooth structure that makes it a Lie groupoid.
Let $\phi: \GR \to \GR$ be an isomorphism of groupoids that
restricts to diffeomorphisms $\GR_1 \to \GR_1$ and $\GR_2 \to
\GR_2$. If the induced map $\phi_*$ map defines a smooth
isomorphism $A(\GR) \simeq A(\GR)$, then $\phi$ is a
diffeomorphism itself.
\end{corollary}

\begin{proof}\
Consider on $\GR$ the smooth structure coming from $\phi$ and
denote by $\maG_\phi$ the resulting Lie groupoid. Then $\phi_*$
establishes an isomorphism of $A(\maG_\phi)$ with $A(\maG)$. By
the previous theorem, Theorem \ref{thm.int}, the smooth structure
$\maG_\phi$ and the given smooth structure on $\maG$ are the same.
Thus $\phi$ is differentiable. Similarly, $\phi^{-1}$ is
differentiable.
\end{proof}

\section{The analytic index\label{a.index}}

We shall need to consider a special class of pseudodifferential
operators on a manifold with poly-cylindrical ends $M_0$ with
compactification $M$, very closely related to the $b$-calculus of
Melrose. For us, it will be convenient to introduce this calculus
using groupoids. For simplicity, {\em we shall assume from
now on that all our manifolds with corners have embedded faces.}

\subsection{Pseudodifferential operators on groupoids}
If $\GR$ is a Lie groupoid with units $M$, then there is
associated to it a pseudodifferential calculus (or algebra of
pseudodifferential operators) $\Psi^{\infty}(\GR)$, whose
operators of order $m$ form a linear space denoted
$\Psi^{m}(\GR)$, $m \in \RR$, such that $\Psi^m(\GR)
\Psi^{m'}(\GR) \subset \Psi^{m + m'}(\GR)$. See \cite{bm-fp, NWX}.
We shall need this construction only for Hausdorff groupoids, so
we assume that $\maG$ is Hausdorff from now on. This calculus is
defined as follows:\ $\Psi^{m}(\GR)$, $m \in \Z$ consists of
smooth families of classical, order $m$ pseudodifferential
operators $(P_x)$, $x \in M$, that are right invariant with
respect to multiplication by elements of $\GR$ and are ``uniformly
supported.'' To define what uniformly supported means, let us
observe that the right invariance of the operators $P_x$ implies
that their distribution kernels $K_{P_x}$ descend to a
distribution $k_P \in I^m(\GR, M)$ \cite{NWX, bm-these}. Then the
family $P = (P_x)$ is called {\em uniformly supported} if, by
definition, $k_P$ has compact support in $\GR$. The right
invariance condition means, for $P = (P_x) \in \Psi^\infty(\GR)$,
that right multiplication $\GR_x \ni g' \mapsto g'g \in \GR_y$
maps $P_y$ to $P_x$, whenever $d(g) = y$ and $r(g) = x$.

We then have the following result \cite{LMN1, bm1, NWX}.

\begin{theorem}\label{thm.quant}\ Let $\GR$ be a Lie
groupoid with units $M$ and Lie algebroid $A = A(\GR)$. The space
$\Psi^{\infty}(\GR)$ is an algebra of pseudodifferential operators
so that there exist surjective principal symbol maps
$\sigma_\GR^{(m)}$ with kernel $\Psi^{m-1}(\GR)$,
\begin{equation*}
    \sigma_\GR^{(m)} : \Psi^m(\GR) \to
    S^m_{cl}(A^*)/S^{m-1}_{cl}(A^*).
\end{equation*}
Also, the algebra $\Psi^{\infty}(\GR)$ acts on $\CI(M)$ such that
$\Psi^{\infty}(\GR)\CIc(M_0) \subset \CIc(M_0)$.
\end{theorem}

In order to define the algebra of pseudodifferential operators
$\Psib$ on our manifold \mc\ $M$ (and acting on $\CIc(M_0)$, where
$M_0$ is our manifold with poly-cylindrical ends), we shall first
define a Lie groupoid $\GR$ canonically associated to $M$, and
then consider the associated pseudodifferential calculus
$\Psi^\infty(\GR)$ of $\GR$ (\cite{LN1, bm1, bm-jfa, NWX}).

Let us denote the hyperfaces of $M$ by $H_i$, $i = 1, \ldots, N$,
as above. Since $M$ has embedded faces, any hyperface $H_i$
has a defining function $\rho_i$, which we shall fix from now on.
To $M$ and $Y = \{\rho_i\}$ we associate the groupoid
\begin{equation}\label{eq.def.GM}
    \tilde G(M; Y) := \{(x,y,\lambda_1,\ldots, \lambda_N)
    \in M \times M \times \rep{N},\
    \rho_i(x) = \lambda_i \rho_i(y),
    \text{ for all } i\}.
\end{equation}

\begin{definition}\label{def.GM}\
We define $G(M)$ to be the $d$-connected component of $\tilde
G(M)$.
\end{definition}

This definition is not canonical, since it depends on the choice
of the defining functions (see \cite{bm-these, bm1} for a canonical
definition). Then $G(M)$ is a Lie groupoid with units $M$. The
operations are the ones induced from the pair groupoid on the
first components and from the group structure on the last
components, as follows:
\begin{equation*}
\begin{gathered}
    d(x, y, \lambda_j) = y, \quad r(x, y, \lambda_j) = x, \;
    \text{ and }\\ (x, y, \lambda_j)(y, z, \lambda_j')
    = (x, z, \lambda_j
\lambda_j').
\end{gathered}
\end{equation*}

We can also consider more general systems $Y$ of functions
$\{\rho\}$ with the property that each function has a
non-degenerate set of zeroes that is a {\em disjoint union} of
hyperfaces of $M$. In particular, we have the following lemma.

\begin{lemma}\label{lemma.open}\
Let $Y = \{\rho_i\}$ and let $Y' = \{\rho_1\rho_2, \rho_3, \ldots
\}$, where the zero sets of $\rho_1$ and $\rho_2$ are disjoint.
Then $\tilde G(M; Y)$ identifies, as a Lie groupoid, with an open
subset of $\tilde G(M; Y')$.
\end{lemma}

\begin{proof}\
The identification of $\tilde G(M; Y)$ with an open subset of
$\tilde G(M; Y')$ is provided by $(x, y, \lambda_i) \to (x, y,
\lambda_1 \lambda_2, \lambda_3, \ldots, )$. Let $Z_j$ be the zero
set of $\rho_j$. The difference between $\tilde G(M; Y)$ and
$\tilde G(M; Y')$ is that $d^{-1}(Z_1 \cup Z_2)$ is $\simeq (Z_1
\cup Z_2)^2 \times \RR$ in the larger groupoid, whereas
$d^{-1}(Z_1 \cup Z_2) \cap \tilde G(M; Y) \simeq (Z_1^2 \cup
Z_2^2) \times \RR$.
\end{proof}

We shall write $C^*(M) := C^*(G(M))$, for simplicity. This
lemma leads right away to the following corollary.

\begin{corollary}\label{cor.open}\
Let $\Omega \subset M$ be an open subset. Fix a system of defining
functions $Y$ for $M$, then $\tilde G(\Omega; Y) \subset \tilde
G(M; Y)$ as an open subset. Consequently, $G(\Omega)$ identifies
canonically with an open subset of $G(M)$ and hence we have a
natural inclusion $\cc{\Omega} \subset \cc{M} = \cc{G(M)}$.
\end{corollary}

\begin{proof}\
Equation \eqref{eq.def.GM} gives
\begin{multline*}
    \tilde G(\Omega; Y) := \{(x,y,\lambda_1,\ldots, \lambda_N)
    \in \Omega \times \Omega \times \rep{N}, \
    \rho_i(x) = \lambda_i \rho_i(y), \text{ for all } i\} \\
    = d^{-1}(\Omega) \cap r^{-1}(\Omega),
\end{multline*}
so $\tilde G(\Omega; Y)$ is indeed an open subset of $\tilde G(M;
Y)$. Passing to the $d$-connected component preserves open
inclusion.
\end{proof}

\begin{remark}\ We clearly have $G(\Omega) \subset d^{-1}(\Omega)
\cap r^{-1}(\Omega)$, but we do not have equality unless there
exists a bijection between the faces of $\Omega$ and those of
$M$.
\end{remark}

We now introduce the class of pseudodifferential operators we are
interested in.

\begin{definition}\label{def.Psib}\ Let $M$ be a manifold
with embedded faces (by our earlier assumption) and $G(M)$ be the
Lie groupoid introduced in Definition \ref{def.GM}, then we let
\begin{equation*}
    \Psib := \Psi^\infty(G(M)).
\end{equation*}
\end{definition}

This algebra is slightly smaller than the one constructed by Melrose
\cite{MelroseScattering}, but for our purposes, it is as good. In
fact, our algebra is the subalgebra of properly supported
pseudodifferential operators in the Melrose's algebra.

We also have that a differential operator $D$ on $M_0$ is in
$\Psib$ \iff it is in $\Diffb$. We also obtain a principal symbol
map
\begin{equation*}
    \sigma_b^{(m)} : \Psi_b^m(M) \to
    S^m_{cl}(A_M^*)/S^{m-1}_{cl}(A_M^*)
\end{equation*}
This new definition of the principal symbol recovers the usual
principal symbol in the interior of $M$, but it also provides
additional information at the boundary. Indeed, the {\em usual}
principal symbol of a differential operator $P \in \Diffb$ is
never invertible at the boundary, so this differential operator
can never be elliptic in the usual sense (unless, of course, it is
actually a function). On the other hand, there are many
differential operators $P \in \Diffb$ whose principal symbol is
invertible on $A_M^*$. An example is provided by $x\pa_x$ on the
interval $[0, \infty)$. This operator has $x\xi$ as its usual
principal symbol, where $(x, \xi) \in T^*[0, \infty) = [0, \infty)
\times \RR$, but in the calculus of $\Psib$, it has principal
symbol $\sigma^{(1)}(x \pa_x) = \xi$.

\subsection{The adiabatic and tangent groupoids}
For the definition and study of the analytic index, we shall need
the adiabatic and tangent groupoids associated to a differentiable
groupoid $\GR$. We now recall their definition and establish a few
elementary properties.

Let $\GR$ be a Lie groupoid with space of units $M$. We construct
both the \textit{adiabatic groupoid } $\adi{\GR}$ and the
\textit{tangent groupoid} $\tgt{\GR}$ (\cite{ConnesNCG, bm-fp, Landsman,
LMN1, Ramazan}). The space of units of $\adi{\GR}$ is $M \times
[0,\infty)$ and the tangent groupoid $\tgt{\GR}$ will be defined as
the restriction of $\adi{\GR}$ to $M \times [0,1]$.

The underlying set of the groupoid $\adi{\GR}$ is the disjoint
union:
\begin{equation*}
        \adi{\GR} = A(\GR) \times \{ 0 \}\, \cup\,
        \GR \times (0,\infty).
\end{equation*}
We endow $A(\GR) \times \{ 0 \}$ with the structure of commutative
bundle of Lie groups induced by its vector bundle structure. We
endow $\GR \times (0,\infty)$ with the product (or pointwise)
groupoid structure. Then the groupoid operations of $\adi{\GR}$ are
such that $A(\GR) \times \{ 0 \}$ and $\GR \times (0,\infty)$ are
subgroupoids with the induced structure. Now let us endow
$\adi{\GR}$ with a differentiable structure. The differentiable
structure on $\adi{\GR}$ is such that
\begin{equation}\label{eq.smooth.gad}
    \Gamma(A(\adi{\GR})) = t \Gamma( A(\GR \times [0, \infty))).
\end{equation}
More precisely, consider the product groupoid $\GR \times [0,
\infty)$ with pointwise operations. Then a section $X \in
\Gamma(A(\GR \times [0, \infty)))$ can be identified with a smooth
function $[0, \infty) \ni t \to X(t) \in \Gamma(A(\GR))$. We thus
require $\Gamma(A(\adi{\GR})) = \{tX(t)\}$, with $X \in
\Gamma(A(\GR \times [0, \infty)))$.

Specifying the Lie algebroid of $A(\adi{\GR})$ completely
determines its differentiable structure \cite{nistorINT}. For
clarity, let us include also an explicit description of this
differentiable structure. Let us consider an atlas
$(\Omega_\alpha)$, consisting of domains of coordinate charts
$\Omega_\alpha \subset \adi{\GR}$ and diffeomorphisms $\phi_\alpha
: \Omega_\alpha \to U_\alpha$, where $U_\alpha$ is an open subset
of a Euclidean space.

Let $\Omega = \Omega_\alpha$ be a chart of $\GR$, such that
$\Omega \cap \GR^{(0)} \neq \emptyset$; one can assume without
loss of generality that $\Omega \simeq T \times U$ with respect to
$s$, and $\Omega \simeq T' \times U$ with respect to $r$. Let us
denote by $\phi$ and $\psi$ these diffeomorphisms.  Thus, if $x
\in U$, $\GR_x \simeq T$, and $A(\GR)_U \simeq \RR^k \times U$.
Let $(\Theta_x)_{x \in U}$ (respectively\ $(\Theta'_x)_{x \in U}$)
be a smooth family of diffeomorphisms from $\RR^k$ to $T$
(respectively\ $T'$) such that $\iota(x)= \phi(\Theta_x(0),x)$
(respectively $\iota(x)=\psi(\Theta'_x(0),x)$, where $\iota$
denotes the inclusion of $\GR^{(0)}$ into $\GR$).

Then $\overline{\Omega}=A(\GR)_U \times \{0\} \cup \Omega \times
(0,\infty)$ is an open subset of $\adi{\GR}$, homeomorphic to
$\RR^k\times U \times \RR_+$ with respect to $s$ and to $r$ as
follows:
\begin{equation*}
\begin{array}{lcll}
\overline{\phi}(\xi,u,\alpha)&=& \left\{\begin{array}{l}
(\phi(\Theta_u(\alpha\xi),u),\alpha)\\ (\xi,u,0)
\end{array}\right.
& \begin{array}{l} \mathrm {\ if\ }\alpha\neq 0\\ \mathrm{\ if\
}\alpha = 0
\end{array} \\
\overline{\psi}(\xi,u,\alpha)&=& \left\{\begin{array}{l} \big(
(\phi(\Theta_u(\alpha\xi),u))^{-1},\alpha\big) \\ (\xi,u,0)
\end{array}\right.
& \begin{array}{l} \mathrm {\ if\ }\alpha\neq 0\\
 \mathrm{\ if\ }\alpha = 0
\end{array}
\end{array}
\end{equation*}
This defines an atlas of $\adi{\GR}$, endowing it with a \cfg\
structure.

Recall that a subset $S$ of the units of a groupoid $\maH$ is
called {\em invariant} \iff\ $d^{-1}(S) = r^{-1}(S)$. Then we
shall denote by $\maH_S := d^{-1}(S) = r^{-1}(S)$ and call it {\em
the restriction of $\maH$ to $S$}. Then $\maH_S$ is also a
groupoid precisely because $S$ is invariant. We shall sometimes
write $\maH\vert_S$ instead of $\maH_S$. For instance, $\tgt{\GR}$,
{\em the tangent groupoid} of $\GR$ is defined to be the
restriction of $\adi{\GR}$ to $M \times [0,1]$.

We have the following simple properties.

\begin{lemma}\label{lemma.times.R}\
Let $\maG$ be a Lie groupoid and $\maH = \maG \times \R^n$ be the
product of $\maG$ with the Lie group $\R^n$ with the induced
product structure. Then $\maH_{ad} \simeq \maG_{ad} \times \RR^n$
\end{lemma}

\begin{proof}\
Let us denote by $M$ the set of units of $\maG$. We have that
$A(\maH) = A(\maG) \oplus \RR^n$ as vector bundles, where the
right copy of $\RR^n$ stands for the trivial bundle with fiber
$\RR^n$ over $M$. Note that, as a set, $A(\maH) = A(\maG) \times
\RR^n$, so our notation will not lead to any confusion.

By definition,
\begin{multline*}
    \maH_{ad} = A(\maH) \times \{ 0 \}\, \cup\,
    \maH \times (0,\infty) = A(\maG) \times \RR^n \times \{ 0 \}\,
    \cup\, \maG \times \RR^n \times (0, \infty) = \\
    \big[ A(\maG) \times \{ 0 \}\,
    \cup\, \maG  \times (0, \infty) \big] \times \RR^n =
    \maG_{ad} \times \RR^n,
\end{multline*}
So we can identify the underlying sets of $\maH_{ad}$ and
$\maG_{ad} \times \RR^n$. However, this identification is not the
bijection we are looking for, because it does not preserve the
differentiable structure. Instead, we define $\phi : \maH_{ad} \to
\maG_{ad} \times \RR^n$ as follows. Let $g \in \maG_{ad}$ and $\xi
\in \RR^n$, so that $(g, \xi) \in \maH_{ad}$, using the previous
identification. Then $\phi(g, \xi) = (g, \xi)$ if $g \in A(\maH)
\times \{ 0 \}$ and $\phi(g, \xi) = (g, t^{-1}\xi)$ if $(g, \xi)
\in \maG \times \{t\} \times \RR^n$.

It remains to check that $\phi$ is differentiable with
differentiable inverse. For this, it is enough to check that it
induces an isomorphism at the level of Lie algebroids, because the
smooth structures on both $\maH_{ad}$ and $\maG_{ad} \times \RR^n$
are defined by their Lie algebroids (see Theorem \ref{thm.int} and
Corollary \ref{cor.int}). Indeed, Equation \eqref{eq.smooth.gad}
gives
\begin{multline*}
    \Gamma(A(\maH_{ad})) = t \Gamma( A(\maH \times [0, \infty)))
    = t \big[ \Gamma( A(\maG \times [0, \infty))) \oplus
    \CI(M \times [0, \infty))^n \big] \\
    \stackrel{\phi_*}{\longrightarrow} t \big[\Gamma( A(\maG
    \times [0, \infty))) \big] \oplus \CI(M \times [0, \infty))^n
    = \Gamma(A(\maG_{ad} \times \RR^n))
\end{multline*}
where $\phi_*$ is an isomorphism.
\end{proof}

This gives the following corollary.

\begin{corollary}\label{cor.prod}\
Let $\maH = \maG \times \RR^n$, as above. We have that
$C^*(\maH_{ad}) \simeq C^*(\maG_{ad}) \otimes \CO(\RR^n)$ and that
$C^*(\maH^{t}) \simeq C^*(\maG^{t}) \otimes \CO(\RR^n)$, the
tensor product being the (complete, maximal) $C^*$--tensor
product.
\end{corollary}

\begin{proof}\
This follows right away from Lemma \ref{lemma.times.R} and the
relation
\begin{equation*}
    C^*(\maG' \times \RR^n) \simeq C^*(\maG')
    \otimes \CO(\RR^n)
\end{equation*}
valid for any locally compact groupoid $\maG'$.
\end{proof}

A similar argument using again Theorem \ref{thm.int} and Corollary
\ref{cor.int} yields the following result.

\begin{lemma}\label{lemma.ad.res}\
Let $\maG$ with a Lie groupoid with units $M$. Let $N \subset M$
be an invariant subset. Assume that $N \subset F$ is an embedded
submanifold of a face $F$ of $M$. Then the restriction operations
and the formation of adiabatic and tangent groupoids commute, in
the sense that we have
\begin{equation*}
    (\maG\vert_N)_{ad} \simeq \maG_{ad}
    \vert_{N \times [0, \infty)}
    \quad \text{and} \quad (\maG\vert_{N^c})_{ad}
    \simeq \maG_{ad} \vert_{N^c \times [0, \infty)}.
\end{equation*}
A similar results holds for the tangent groupoids.
\end{lemma}

\begin{proof}\ Again, it follows from the definition that
$(\maG\vert_N)_{ad}$ and $\maG_{ad}\vert_{N \times [0, \infty)}$
have the same underlying set. Moreover, this canonical bijection
is a groupoid isomorphism. From the definition of the Lie
algebroid of the adiabatic groupoid, if follows that this
canonical bijection is also differentiable with differentiable
inverse, because it induces an isomorphism of the spaces of
sections of the corresponding Lie algebroids. See Theorem
\ref{thm.int} and, especially, Corollary \ref{cor.int}.
\end{proof}

We then obtain the following corollary.

\begin{corollary}\label{cor.ad.res}\ With the notations of the
above lemma we have a short exact sequence
\begin{equation*}
    0 \to C^*(\tgt{\maG\vert_{N^c}}) \to C^*(\tgt{\maG}) \to
    C^*(\tgt{\maG\vert_{N}}) \to 0,
\end{equation*}
and a similar exact sequence for the adiabatic groupoid.
\end{corollary}

\subsection{The analytic index}
We now give two definitions of the analytic index. The first
definition is based on a generalization of Connes' tangent
groupoid \cite{ConnesNCG}. The second definition is based on the
boundary map of the six-term exact sequence in \kth\ induced by
the symbol map. Both definitions will be needed in what follows.

For each $t\in [0,1]$, $M \times \{t\}$ is a closed invariant
subset of $M \times [0, \infty)$, and hence we obtain an {\em
evaluation map}
\begin{equation*}
    e_t : \cc{\tgt{\GR}} \to \cc{\tgt{\GR}_{M \times \{t\}}}.
\end{equation*}
By abuse of notation, we shall sometimes denote also by $e_t$ the
induced map in \kth.

Let us also notice that the decomposition
\begin{equation*}
    M \times [0,1] = M \times \{0 \} \cup M \times (0,1]
\end{equation*}
into two closed, invariant subspaces gives rise to an exact
sequence
\begin{equation}\label{eq.exact1}
    0 \to \cc{\tgt{\GR}_{M \times (0,1]}} \to \cc{\tgt{\GR}}
     \xrightarrow{e_0} C^*(A(\GR)) \to 0,
\end{equation}
This leads to the following six-terms exact sequence
in \kth:
\begin{equation*}
  \begin{CD}
    K_0( \cc{\tgt{\GR}_{M \times (0,1]}})  @>>> K_0( \cc{\tgt{\GR}})
    @>>> K_0(C^*(A(\GR)))\\ @AAA @. @VVV\\
    K_1(C^*(A(\GR)) )  @<<< K_1( \cc{\tgt{\GR}}) @<<<
    K_1(\cc{\tgt{\GR}_{M \times (0,1]}}).\\
  \end{CD}
\end{equation*}

We have $\tgt{\GR}_{M \times (0,1]}=\GR \times (0,1]$ and hence
$\cc{\tgt{\GR}_{M \times (0,1]}} \simeq \cc{\GR}\otimes \CO((0,1])$.
In particular, $K_*( \cc{\tgt{\GR}_{M \times
(0,1]}})=K_*(\cc{\GR}\otimes \CO((0,1])=0$. Thus the evaluation
map $e_0$ is an isomorphism in \kth.

The $C^*$-algebra $C^*(A(\GR))$ is commutative and we have
$C^*(A(\GR)) \simeq \CO(A^*(\GR))$. Therefore
$K_*(C^*(A(\GR)))=K^*(A^*(\GR))$. In turn, this isomorphism allows
us to define the \textit{analytic index} $\ind_a$ as the
composition map
\begin{equation}\label{def.an.index1}
    \ind_a^\GR = e_1 \circ e_0^{-1} : K^*(A^*(\GR))
    \to K_*(\cc{\GR}),
\end{equation}
where $e_1 : \cc{\tgt{\GR}} \to \cc{\tgt{\GR}_{M \times \{1\}}} =
\cc{\GR}$ is defined by the restriction map to $M \times \{1\}$.

The definition of the analytic index gives the following.

\begin{proposition}\label{prop.comp}\
Let $\GR$ be a Lie groupoid with Lie algebroid $\pi : A(\GR) \to
M$. Also, let $N \subset F \subset M$ be a closed, invariant
subset which is an embedded submanifold of a face $F$ of $M$. Then
the analytic index defines a morphism of the six-term exact
sequences associated to the pair $(A^*(\maG), \pi^{-1}(N))$ and to
the ideal $\cc{\maG_{N^c}} \subset \cc{\maG}$, $N^c := M
\smallsetminus N$
\begin{equation*}
  \begin{CD}
    K^0(\pi^{-1}(N^c))  @>>> K^0(A^*(\maG))
    @>>> K^0(\pi^{-1}(N)) @>>> K^1(\pi^{-1}(N^c))\\
    @VVV @VVV @VVV @VVV\\
    K^0(\cc{\maG_{N^c}})  @>>> K^0(\cc{\maG})
    @>>> K^0(\cc{\maG_{N}}) @>>> K^1(\cc{\maG_{N^c}})\\
  \end{CD}
\end{equation*}
\end{proposition}

\begin{proof}\
The six-term, periodic long exact sequence in \kth\ associated to
the pair $(A^*(\maG), \pi^{-1}(N))$ is naturally and canonically
isomorphic to the six-term exact sequence in \kth\ associated to
the pair $\CO(A^*_{M \smallsetminus N}) \subset \CO(A^*(\maG))$
consisting of an algebra and an ideal in that algebra. Since $e_0$
and $e_1$ also induce morphisms of pairs (algebra, ideal), the
result follows from Corollary \ref{cor.ad.res}, the naturality of
the six-term exact sequence in \kth, and the definition of the
analytic index \eqref{def.an.indexM}.
\end{proof}

For $M$ a smooth manifold with embedded faces, we have
$A(G(M))=A_M$. Recall that $\cc M := C^*(G(M))$. Then the analytic
index becomes the desired map
\begin{equation}\label{def.an.indexM}
    \ind_a^M : K^*(A^*_M) \to K_*(\cc{M}).
\end{equation}

\begin{remark}\ Assume $M$ has no corners (or boundary).
Then $G(M) = M \times M$ is the product groupoid and hence
$\Psi^\infty(G(M)) = \Psi^\infty(M)$. In particular, $C^*(M) :=
C^*(G(M)) \simeq \maK$, the algebra of compact operators on $M$.
In this case $K_0(C^*(M))=\ZZ$, and $\ind_a$ is precisely the
analytic index as introduced by \cite{AS1}. This construction
holds also for the case when $M$ is not compact, but we have to
use pseudodifferential operators of order zero that are
``multiplication at infinity,'' as in \cite{carvalho}.
\end{remark}

We now turn to the application mentionned in the introduction
to obstructions to finding invertible perturbations by
regularizing operators. Let $\mathfrak A(\GR)$ be the envelopping
$C^*$--algebra of $\Psi^0(\GR)$ and let $S^*\GR \subset A^*(\GR)$
be the subset of vectors of length one. The second method is based
on the exact sequence
\begin{equation}\label{eq.exact2}
    0 \to C^*(\GR) \to \mathfrak A(\GR)
     \xrightarrow{\sigma^{(0)}} C(S^*\GR) \to 0.
\end{equation}
The six-term exact sequence in $K$--theory associated to this
exact sequence yields a boundary map
\begin{equation}\label{def.an.index2}
    \ind_a^\GR := \pa : K^{*+1}(S^*\GR) \to K_*(\cc \GR).
\end{equation}

Considering the exact sequence
$$0 \to C_0(A(\GR)) \to C(B(\GR)) \to C(S(\GR))\to 0,$$
where $B(\GR)$ is the ball bundle, we get a map $b:K^{*+1}(S^*\GR)\to
K^*(A^*M)$, and $\ind_a^M\circ b= \ind_a^\GR $.

This second definition of the analytic index has the advantage
that it leads to the following theorem. Let $\Psib$ be as in
Definition \ref{def.Psib} and
\begin{equation*}
    \Psi_b^m(M; E_0, E_1) = e_1M_N(\Psi_b^m(M))e_0,
\end{equation*}
where $e_j$ is the orthogonal projection onto the subbundle $E_j
\subset M \times \CC^N$, $j = 1, 2$.

\begin{theorem}\label{thm.perturb.M}\
Assume $M$ is a connected manifold with corners such that all its
faces have positive dimension. Let $P \in \Psi_b^m(M; E_0, E_1)$
be an elliptic pseudodifferential operator acting between sections
of two vector bundles $E_0, E_1 \to M$. Then there exists $R \in
\Psi_b^m(M; E_0, E_1)$ such that $P + R$ is invertible if, and
only if, $\ind_a^M(\sigma_P) = 0$, where $\sigma_P$ is the
principal symbol of $P$.
\end{theorem}

\begin{proof}\ If $E_0 = E_1$ is the trivial bundle, then the
proof is the same as that of Theorem 4.10 of \cite{NT1}.

Let us observe as in \cite{Moroianu} that $E_0 \simeq E_1$ because
$TM$ has a non-zero section. By embedding $E_0$ into a trivial
bundle, we can therefore assume that $P \in M_N(\Psi_b^m(M))$. The
index $\ind_a^M (\sigma_P) = \pa [\sigma_P] \in K_0(C^*(M))$ is
therefore defined as in Equation \eqref{def.an.index2}. Moreover,
this definition is independent of the isomorphism $E_0 \simeq E_1$
and of the embedding of $E_0$ into a trivial bundle. This reduces
the proof to the case of a trivial bundle.
\end{proof}

\section{Properties of the analytic
index\label{sec.prop}}

We now prove some  results on the analytic index whose definition
was recalled in the previous section. In this and the following
sections, we continue to assume that our manifolds with corners
have embedded faces. Recall that $G(M)$ is the Lie groupoid
associated to a manifold with corners $M$ in Definition
\ref{def.GM}.

Also, recall that we have denoted $A_M := A(G(M))$ and $C^*(M) :=
C^*(G(M))$. In particular, $\Gamma(A_M) = \maV_b$, the space of all
vector fields tangent to the faces of $M$.

\subsection{The role of faces}
Let $F \subset M$ be a face of a manifold with embedded faces $M$.
Recall that in our terminology, ``face'' always means ``closed
face.'' Then $F$ is a closed, invariant subset of $M$
(``invariant'' here is with respect to the action of $G(M)$ on its
units).

\begin{lemma}\label{lemma.face}\
Let $\pi : A_M \to M$ be the Lie algebroid of $G(M)$. Then, for
any face $F \subset M$ of codimension $k$, we have isomorphisms
\begin{equation*}
    \pi^{-1}(F) \simeq A_F \times \RR^k \quad \text{and}
    \quad G(M)\vert_{F} \simeq G(F) \times \RR^k.
\end{equation*}
\end{lemma}

\begin{proof}\
We have that $\pi^{-1}(F) \simeq A_F \oplus \pi^{-1}(F)/A_F$. The
choice of defining functions for the $k$ hyperfaces containing $F$
then gives an isomorphism $\pi^{-1}(F)/A_F \simeq F \times \RR^k$.
The last part follows from the definitions of $G(M)$ and $G(F)$,
for the latter using the defining functions of $M$ that are
non-zero on $F$.
\end{proof}

For the simplicity of the notation, let us denote by $\maG_1 :=
G(M) \vert_{F} \simeq G(F) \times \RR^k$. The analytic indices
defined in Equations \eqref{def.an.index1} and
\eqref{def.an.indexM} for this restriction groupoid are then
identified by the following lemma.

\begin{lemma}\label{lemma.a.i.f}\
For any face $F \subset M$ of codimension $k$ we have a
commutative diagram
\begin{equation*}
\begin{CD}
    K^{*+k}(A^*(\maG_1)) @>\ind_a^{\maG_1}>>
    K_{*+k}(C^*(\maG_1))\\ @AAA @AAA \\
    K^*(A^*_F) @>\ind_a^{F}>> \kc {F},\\
\end{CD}
\end{equation*}
where the vertical arrows are the periodicity isomorphisms.
\end{lemma}

\begin{proof}\
Recall first that the periodicity isomorphism $K_{*}(\mfk A)
\simeq K_{*+k}(\mfk A \otimes \CO(\RR^k))$ is natural in $\mfk A$,
for any $C^*$-algebra $\mfk A$. The result then follows from this
observation combined with Lemma \ref{lemma.face} above and
Corollary \ref{cor.prod}.
\end{proof}

The following proposition is one of the main steps in the proof of
our topological index theorem, \ref{thm.topological}.

\begin{proposition}\label{prop.2}\
Let $X$ be a \mec\ such that each open face of $X$ is
diffeomorphic to a Euclidean space. Then the analytic index
\begin{equation*}
    \ind_a^X : K^*(A^*_X) \to K_*(\cc{X})
\end{equation*} (defined in
Equation \eqref{def.an.indexM}) is an isomorphism.
\end{proposition}

\begin{proof}\
We shall proceed by induction on the number of faces of $X$. Let
$G(X)$ be the groupoid of $X$, as before. Let $F \subset X$ be a
face of minimal dimension. Then $F$ is an invariant subset of $X$.
Denote as above by $\maG_1$ the restriction of $G(X)$ to $F$. In
particular, $F$ will be a smooth manifold {\em without corners}.

Note that our assumptions imply that $F \simeq \RR^{n-k}$, where
$n$ is the dimension of $M$ and $k$ is the codimension of $F$.
Therefore the analytic index $\ind_a^F$ is an isomorphism. In
particular, our result is valid if $X$ has exactly one face. Lemma
\ref{lemma.a.i.f} then shows that $\ind_a^{\maG_1}$ is an
isomorphism as well.

We shall complete the proof using Proposition \ref{prop.comp} with
$N = F$ and $X = M$ as follows. Let us consider the six-term,
exact sequence in \kth\ associated to the pair $(A_X^*,
A^*(\maG_1))$, where $A^*(\maG_1) = \pi^{-1}(F)$, with $\pi :
A^*_X \to X$ the canonical projection. Also, let us consider the
six-term, exact sequence in \kth\ associated to the pair $\cc{F^c}
\subset \cc{X}$, $F^c = X \smallsetminus F$. Proposition
\ref{prop.comp} states that the analytic index defines a morphism
of these two exact sequences. The analytic index for the quotient
(\ie\ $\ind_a^{\maG_1}$) has just been proved to be an
isomorphism. The analytic index for $F^c$ is an isomorphism by the
induction hypothesis. The Five Lemma then shows that the remaining
two analytic index morphisms are also isomorphisms. The proof is
now complete.
\end{proof}

\begin{remark}
The above proposition can be regarded as a Baum--Connes
isomorphism for manifolds with corners.
\end{remark}

Another important step in the proof of our topological index
theorem for manifolds with corners is the proof of the following
proposition. We first need a lemma.

\begin{lemma}\label{lemma.prop.3}\
Let $\iota : M \to X$ be an embedding of manifolds with corners.
Then $\iota$ defines a natural morphism
\begin{equation}
    \iota_K : \kc{M} \to  \kc{X}.
\end{equation}
\end{lemma}

\begin{proof}\ The morphism $\iota_K$ will be defined by a
Kasparov $\cc{M}-\cc{X}$ bimodule $\maH$. Let us first identify
$M$ with a closed submanifold with corners of $X$. Let $X'
\subset X$ be a small neighborhood of $M$ in $X'$ such that every
face of $X'$ intersects $M$ and the map $F \to F \cap M$
establishes a bijection between the faces of $X'$ and those of
$M$. For instance, $X'$ could be a tubular neighborhood of $M$ in
$X$. Since $X' \subset X$ is open, there exists a natural
inclusion $\cc{X'} \to \cc{X}$ of $C^*$-algebras, by Corollary
\ref{cor.open}. It is enough then to construct a natural morphism
$\kc{M} \to \kc{X'}$. We can therefore replace $X$ with $X'$ and
assume that there exists a bijection between the faces of $M$ and
those of $X$.

Let us introduce
\begin{equation*}
    \Omega := r^{-1}(M) \subset G(X),
\end{equation*}
that is, $\Omega$ is the set of elements of $G(X)$ whose range is
in $M$. We endow $\Omega$ with the two maps $r : \Omega \to M$ and
$d: \Omega \to X$ induced by the range and domain maps of the \gr\
$G(X)$. These are continuous and open since $G(X)$ is a \cfg. The
groupoid $G(X)$ acts on the right on $\Omega$ by multiplication.
We shall define $\maH$ as a Hilbert module completion of
$\CO(\Omega)$ as in \cite{renault-equiv}.

We need next to define the action of $\cc{G(M)}$. This action will
be by compact operators (``compact'' here is used in the sense of
Hilbert modules) and will come from an action of $G(M)$ on
$\Omega$. To define this action, we first prove that
\begin{equation}\label{eq.prop.GM}
    G(M) = r^{-1}(M) \cap d^{-1}(M).
\end{equation}
Indeed, our assumptions imply that for any hyperface $H$ of $X$
with defining function $\rho = \rho_H$, the function
$\rho\vert_{M}$ is a defining function of $M$. (We are, of course,
using the fact that $M$ and $H$ intersect transversely.) Our
assumptions imply in fact more, they imply that every defining
function of $M$ is obtained in this way. We can therefore
establish a bijection between the defining functions $\rho$ of
(the hyperfaces of) $X$ and the defining functions of (the
hyperfaces of) $M$. Let us use these functions in the definition
of $G(M)$, namely in Equation \eqref{eq.def.GM}. This, together
with the fact that there is a bijection between the open faces of
$M$ and the open faces of $X$, proves then Equation
\eqref{eq.prop.GM}.

Equation \eqref{eq.prop.GM} then allows us to define the action of
$G(M)$ on $\Omega$ by left composition. Indeed, if $\gamma \in
G(M)$ and $\gamma' \in \Omega$ with $r(\gamma')=d(\gamma)$, then
$r(\gamma\gamma')=r(\gamma) \in M$ so that $\gamma\gamma' \in
\Omega$. Also $G(X)$  acts on $\Omega$ by right composition: if
$\gamma \in G(X)$ and $\gamma' \in \Omega$ with
$d(\gamma')=r(\gamma)$, then $r(\gamma'\gamma)=r(\gamma') \in M$
so that $\gamma'\gamma \in \Omega$.

Then $\maH$ defines an element in $\Theta \in KK_0(\cc{M},\cc{X})$
(in fact, even an imprimitivity module) and hence the Kasparov
product with $\Theta$ defines the desired morphism $\iota_K :
\kc{M} \to \kc{X}$.
\end{proof}

\begin{remark}\label{rem.prop.3}\ Let us spell out explicitly a
conclusion of the above proof. If $M \subset X$ is an open subset
and $\iota$ denotes the inclusion then $\cc{M} \subset \cc{X}$, by
Corollary \ref{cor.open}, and $\iota_K$ is simply the morphism
associated to this inclusion of $C^*$-algebras.
\end{remark}

\begin{proposition}\label{prop.3}
Let $\iota: M \to X$ be a closed embedding of manifold with
corners. Assume that, for each open face $F$ of $X$, the
intersection $F \cap M$ is a non-empty open face of $M$ and that
every open face of $M$ is obtained in this way. Then $\iota_K :
\kc{M} \to  \kc{X}$ is an isomorphism.
\end{proposition}

\begin{proof} Recall from \cite{renault-equiv}
that two locally compact groupoids $G$ and $H$ are equivalent
provided there exists a topological space $\Omega$ and two
continuous, surjective open maps $r : \Omega \to \g0$ and $d :
\Omega \to \h0$ together with a left (respectively right) action
of $G$ (respectively $H$) on $\Omega$ with respect to $r$
(respectively $d$), such that $r$ (respectively $d$) is a
principal fibration of structural \gr\ $H$ (respectively $G$).

An important theorem of Muhly--Renault--Williams states that if
$G$ and $H$ are equivalent, then $\kc {G} \simeq \kc {H}$
\cite{MRW}. More precisely, $\CO(\Omega)$ has a completion to an
Hilbert module that establishes a strong Morita equivalence
between $\cc{G}$ and $\cc{H}$ (this is the imprimitivity
module defining $\Theta$ in the proof of Lemma
\ref{lemma.prop.3}). This strong Morita equivalence is then known
to imply the stated isomorphism $\kc {G} \simeq \kc {H}$
\cite{Rieffel}.

To prove our result, it is therefore enough to show that the space
$\Omega := r^{-1}(M)$ (considered also in the proof of Lemma
\ref{lemma.prop.3}) establishes an equivalence between $G(M)$ and
$G(X)$.

Proving that $r$ is a principal fibration of structural \gr\
$G(X)$ amounts, by definition, to proving that, for any $x$ in
$M$, if $\omega$ and $\omega'$ are in $r^{-1}(x)$ in $\Omega$,
there exists $\gamma \in G(X)$ such that $\omega \gamma =
\omega'$, and that the action of $G(X)$ is free and proper. The
first condition is clear:\ $r(\omega) = r(\omega') = x$ so that
$\gamma = \omega^{-1}\omega' \in G(X)$ exists. Besides, the action
is free. Indeed, if $\omega\gamma=\omega$, then $r(\gamma) =
d(\omega)$ (so they are composable) and $d(\gamma) =
d(\omega\gamma) = d(\omega)$, so that $\gamma$ is a unit. The
action is proper. Indeed, the map
\begin{equation*}
     \phi: G(X) * \Omega \to \Omega\times \Omega,
     \qquad (\gamma,\omega) \mapsto (\omega\gamma,\omega)
\end{equation*}
(where $ G(X) * \Omega$ is the set of composable arrows) is a
homeomorphism onto its image, which is the fibered product $\Omega
\times_M \Omega$ with respect to $r$.

Similarly, let us check that $s$ is a principal fibration with
structural \gr\ $G(M)$. Assume that $d(\omega) = d(\omega') \in M$,
then $\omega' = \gamma \omega$, with $\gamma = 
\omega' \omega^{-1}$. Let us recall from the proof of Lemma
\ref{lemma.prop.3} that $G(M) = d^{-1}(M) \cap r^{-1}(M)$. Hence
$\gamma \in G(M)$. The proof is now complete.
\end{proof}

Let $M$ be a smooth compact manifold (without corners). Then the
inclusion of a point $k : pt \hookrightarrow M$ satisfies the
assumptions of the above Proposition. The imprimitivity module in
this case establishes the isomorphism $C^*(pt) \otimes \maK \simeq
C^*(M)$ and hence also the isomorphisms $K_*(C^*(M)) \simeq
K_*(C^*(pt)) = K_*(\CC) \simeq \ZZ$.

More generally, let $M$ be a compact manifold with corners
and  $\pi : X \to M$ be a smooth fiber bundle. Recall then that
the {\em fiberwise product groupoid} $\GR := X \times_M X$ was
defined in the first section. The Lie algebroid of $\GR$ is
$T_{\pi}X$, the vertical tangent bundle to $\pi : X \to M$. A
simple calculation shows that $C^*(\GR)$ is a continuous field of
$C^*$-algebras over $M$ with fibers compact operators on $L^2$ of
the fibers of $\pi$. Therefore $C^*(\GR)$ is canonically Morita
equivalent to $C(M)$. If $X \to M$ has a cross section (as in the
cases when we shall use this construction in our paper), we can
also obtain this Morita equivalence from the inclusion $M \subset
X \times_M X$ given by this cross-section, which is an equivalence
of groupoids. In any case, we obtain an isomorphism
\begin{equation}\label{theta}
    K_*(C^*(\GR)) \simeq K_*(C(M)) = K^*(M),
    \quad \GR := X \times_M X,
\end{equation}
which we shall often use to identify these groups. For instance,
the analytic index associated to $\GR = X \times_M X$ becomes a
map $\ind_a^{\GR} : K^*(T^*_{\pi}X) \to K^*(M) = K_*(C(M))$. Then
we have the following well known result, whose proof we sketch for
the benefit of the reader.

\begin{proposition}\label{prop.family}
Let $\pi : X \to M$ be a smooth fiber bundle with $X$
compact and $\GR := X \times_M X$ be the fiberwise product
groupoid. Let $[a] \in K^0(T^*_{\pi}X)$ be a $K$-theory class
represented by an endomorphism $a : E \to F$ of vector bundles
over $T^*_{\pi}X$ that are pull-backs of vector bundles on $X$ and
are such that $a$ is an isomorphism outside a compact set. Assume
that $a$ is homogeneous of some positive order and let $P_a$ be a
family of (elliptic) pseudodifferential operators along the fibers
of $\pi$ with principal symbol $a$. Then $\ind_a^{\GR}([a]) \in
K^0(M)$ coincides with the family index of $P_a$.
\end{proposition}

\begin{proof}\
Let $B$ be a $C^*$-algebra. We shall denote by $B^+$ the algebra
with an adjoint unit. For any $b \in B$, we shall denote by $p_b
\in M_2(B^+)$ the graph projection associated to $b$, that is
\begin{eqnarray}
        \label{eqgraph}
        p_b = \left[\begin{array}{cc}1- e^{-B^*B}&\tau(B^*B)B^*\\ 
        \tau(BB^*)B & e^{-BB^*}
        \end{array}\right].
\end{eqnarray}
where  $\tau$ is a
smooth, even function on $\RR$ satisfying
$\tau(x^2)^2x^2=e^{-x^2}(1-e^{-x^2})$.

We can assume that $a$ is polynomial, homogeneous of degree one, given
by the principal symbol of a family $D$ of first order differential
operators. Let $p_a$ be the graph projection associated to $a$, $p_a
\in M_N(\CO(T^*_{\pi}X)^+)$. Then we can extend the family $a$ to the
family $p_{tD} \in M_N(C^*(\tgt{\GR})$. Therefore, the analytic index of
$[a]$ is the class of $p_{D} \in K_0(C^*(\GR))$. We can furthermore
perturb $tD$ with a regularizing family that vanishes in a
neighborhood of $t = 0$ without changing the class of $p_D$.  We can
hence assume that $D$ has a kernel of constant dimension
\cite{AS1}. Then it is known that the class of $p_D$ is $[\ker D] -
[\ker D^*]$, that is, the family index of $D$.
\end{proof}

We now prove a corollary that will be needed in the proof of
Proposition \ref{prop.diag1}. Let $\pi : U \to X$ be a vector
bundle over a compact manifold with corners. Consider the fibered
product groupoid $\GR : = U \times_X U$ as above. Then $T_\pi^* U
= U \oplus U^*$ as vector bundles over $X$, and hence $i_! :
K^*(X) \to K^*(T_\pi^* U)$ is defined.

\begin{corollary}\label{cor.family}\
We have that $\ind_a^\GR \circ i_!$ is the identity map of
$K^*(X)$ (or, more precisely, the inverse of the isomorphism of
Equation \eqref{theta}).
\end{corollary}

\begin{proof}\ 
Both $\ind_a^\GR$ and $i_!$ are $K^*(X)$-linear. The same argument as
in the proof of Proposition \ref{prop.family} (or by compactifying $U$
fiberwise to a sphere bundle and using then Proposition
\ref{prop.family}), we obtain that $\ind_a^\GR \circ i_! (1)$ is the
index of the family of Dirac operators on the fibers of $U$ coupled
with the potential given by Clifford multiplication with the
independent variable. Since the equivariant index of the coupled Dirac
operator is $1$, we obtain that $\ind_a^\GR \circ i_! (1) = 1$ as in
\cite{AS1}. See \cite{hor3} for a simple proof of the facts needed
about the coupled Dirac operator.
\end{proof}

\begin{remark}\label{rem.top}\ Let us take $X = pt$, to be reduced
to a point, and let us identify $U$ with $\RR^N$, for some $N$.
Then the above Corollary states, in particular, that $i_! : \ZZ =
K^*(pt) \to K^*(T\RR^n)$ and $\ind_a^{\RR^n} : K^*(T\RR^n) \to \ZZ
= K^*(pt)$ are inverse to each other.
\end{remark}

\section{Commutativity of the
diagram\label{commutativity-diagram}}

In this section we shall prove a part of our topological index theorem,
Theorem \ref{thm.topological}, involving an embedding $\iota : M
\to X$ of our manifold with corners $M$ into another manifold with
corners $X$. This theorem amounts to the fact that the diagram
\eqref{diag.I} is commutative. In order to prove this, we shall
first consider a tubular neighborhood
\begin{equation}\label{eq.def.jk}
    M \overset{k}{\hookrightarrow} U \overset{j}{\hookrightarrow} X
\end{equation}
of $M$ in $X$. The diagram (\ref{diag.I}) is then decomposed into
the two diagrams below, and hence the proof of the commutativity
of the diagram (\ref{diag.I}) reduces to the proof of the
commutativity of the two diagrams below, whose morphisms are as
follows:\ the morphisms $k_K$ and $j_K$ are defined by Lemma
\ref{lemma.prop.3}, the morphism $k_!$ is the push-forward
morphism, and $j_*$ is the morphism in \kth\ defined by an open
embedding. Recall that $A_M = A(G(M))$ and $C^*(M) := C^*(G(M))$.

\begin{equation}\label{diag2}
\begin{CD}
    \kc{M} @>{k_K}>> \kc{U} @>{j_K}>>\kc{X}\\
    @A{\ind_a^M}AA @A{\ind_a^U}AA @AA{\ind_a^X}A\\ K^*(A^*_M)
    @>{k_{!}}>> K^*(A^*_U) @>{j_*}>> K^*(A^*_X)\\
\end{CD}
\end{equation}

The commutativity of the left diagram is part of the following
proposition.

\begin{proposition}\label{prop.diag1}
Let $\pi : U \to M$ be a vector bundle over a manifold with
corners $M$ and let $k : M \to U$ be the ``zero section''
embedding. Then the following diagram commutes:
\begin{equation}\label{diag3}
\begin{CD}
    \kc{M} @>{k_K}>{\simeq}> \kc{U}\\
    @A{\ind_a^M}AA @A{\ind_a^U}AA\\
    K^*(A^*_M) @>{\simeq}>{k_{!}}> K^*(A^*_U)
\end{CD}
\end{equation}
\end{proposition}

\begin{proof}\
We shall prove this result using a double deformation groupoid $\maG$,
which is a Lie groupoid with units $U \times [0,1]^2$. This
groupoid is such that the projection $U \times [0, 1]^2 \to [0,
1]^2$ extends to a groupoid morphism $\GR \to [0, 1]^2$ the latter
being considered as a space, \ie a groupoid equal to its units. In
other words, if $d$ and $r$ are the domain and range of $\GR$, then
$d(g)$ and $r(g)$ have the same projection in $[0, 1]^2$.

As a set, $\maG := \GR_1 \sqcup \GR_2 \sqcup
\GR_3$ ($\sqcup$ denotes the disjoint union), where 
\begin{equation*}
\begin{gathered}
    \maG_1 := A_U \times \{0\}\times [0, 1], \quad
    \maG_2 := A_M \times_M U \times_M U \times (0,1]\times
    \{0\},\\ \quad \maG_3 := G(U) \times (0,1]\times (0,1].
\end{gathered}
\end{equation*}

\begin{figure}[htbp]
\centering
  \input{mn3.pstex_t}
 \caption{The groupoid $\GR$}
  
\end{figure}

\begin{lemma}
  The groupoid $\GR$ is a \cfg.
\end{lemma}
\begin{proof}
We prove this by applying theorem \ref{thm.int} which states that there exists at most one
differentiable structure on $\GR$ compatible with the groupoid
structure (that is, making it a Lie groupoid). We define a Lie
algebroid structure, and Lie groupoid structures for each of the
subgroupoids $\GR_1, \GR_2, \GR_3$, and last we show how to get a
\cfg\ structure compatible with this data.

Let us first describe $A(\GR)$, the Lie algebroid of the Lie groupoid
that we want to construct. Recall that the Lie algebroid associated to
$U$, $A_U$, or ``compressed tangent bundle,'' is such that
$\Gamma(A_U)$ consists of the smooth vector fields on $U$ that are
tangent to all faces of $U$. Then $\pi : U \to M$ induces a map $\pi_*
: A_U \to A_M$.

Let $A_U = A_v \oplus A_h$ be a decomposition of $A_U$ into vertical
and horizontal components, so that $A_v$ is the kernel of $\pi_* : A_U
\to A_M$ and $A_h \simeq \pi^*(A_M)$. Then we obtain the decomposition
\begin{equation*}
    \Gamma(U \times [0, 1]^2; A_U) = \Gamma(U \times [0, 1]^2; A_v)
    \oplus \Gamma(U \times [0, 1]^2; A_h).
\end{equation*}
Let us regard the above spaces of smooth sections $X$ as families of
sections $X(s, t) \in \Gamma(A_U)$ parameterized by $(s, t) \in [0,
1]^2$. Then $A(\GR)$, the Lie algebroid of the Lie groupoid that we
want to construct is defined, as a set, by
\begin{equation*}
    \Gamma(A(\maG)) = s\Gamma(U \times [0, 1]^2; A_v) \oplus
    st\Gamma(U \times [0, 1]^2; A_h).
\end{equation*}
Then we define the Lie bracket to be the pointwise Lie bracket coming
from the Lie bracket on sections of $\Gamma(A_U)$:\ $[X, Y](s, t) =
[X(s, t), Y(s, t)]$.

Let us now define the Lie groupoid structure of each subgroupoid as follows. 
First,  $\GR_1$ is a commutative Lie groupoid, with operations
defined by the vector bundle structure on $A_U$. It integrates the
restriction of $A(\maG)$ to $\{s = 0\} = \{0\} \times [0, 1]$. The
groupoid $\GR_2$ is the fibered product of $A_M \to M$ and of the
fiberwise pair groupoid $U \times_M U$ considered in Proposition
\ref{prop.family} (in the case at hand, $X = U$). More precisely,
let $\pi_M : A_M \to M$ be the canonical projection, then
\begin{multline*}
    A_M \times_M U \times_M U \times (0,1]\times \{0\} =
    \{(\xi, u_1, u_2, s, 0) \in A_M \times U
    \times U \times (0,1]\times \{0\}, \\
    \pi_M(\xi) = \pi(u_1) = \pi(u_2)\},
\end{multline*}
with the product $(\xi, u_1, u_2, s, 0)(\xi', u_2, u_3, s, 0)
=(\xi + \xi', u_1, u_3, s, 0)$. The factor $(0,1]\times \{0\}$
therefore plays just the role of a space of parameters. The last
groupoid, $\GR_3$ is the product of the groupoid $G(U)$ associated
to the manifold with corners $U$ (Definition \ref{def.GM}) with
the space $(0, 1] \times (0, 1]$, which again plays just the role
of a space of parameters. It integrates the restriction of
$A(\maG)$ to $(0, 1] \times (0, 1]$.

Theorem \ref{thm.int} states that there exists at most one
differentiable structure on $\GR$ compatible with the groupoid
structure (that is, making it a Lie groupoid). The union $\GR_1
\sqcup \GR_3$ has the differentiable structure of the tangent
groupoid of $G(U)$, product with $(0, 1]$. Due to the local
structure of the deformations involved, in order to prove that
$\GR$ does indeed have a smooth structure, it is enough to assume
that $U \to M$ is a trivial bundle, in which case the resulting
$\GR$ is seen to be smooth as follows.

Let $U = M \times \RR^n$. Then $A_U = A_M \times T\RR^n$ and
$\maG_2 = A_M \times (\RR^n \times \RR^n) \times (0,1] \times
\{0\}$, with $(\RR^n \times \RR^n)$ being the pair groupoid.
Deforming just in the fiberwise direction means, in this case,
that we deform the tangent space to $\RR^n$ to its tangent
groupoid. We then define the topology on $\maG_1 \cup \maG_2$,
which is a groupoid with units $\{0\}\times [0, 1] \cup
(0,1]\times \{0\} \supset [0,1]\times \{0\}$, such that the
restriction to $[0,1]\times \{0\}$ is $A_M \times \tgt{G(\RR^n)}$,
with $\tgt{G(\RR^n)}$ denoting the tangent groupoid of the smooth
manifold $\RR^n$. This is possible since the tangent groupoid of
$\RR^n$ at $t = 0$ is $T\RR^n$ and the restriction of $\GR_1$ at
$t = 0$ is $A_M \times T\RR^n$. 

On the other hand, since $U=M\times \RR^n$, the groupoid $G(U)$ is
just the product of $G(M)$ by the pair groupoid $\RR^n \times \RR^n$.
The union $\GR_2 \sqcup \GR_3$ is $\tgt{G(M)} \times (\RR^n \times \RR^n)
\times (0, 1]$, with the factor $(0, 1]$ corresponding to the variable
$s$.
\end{proof}

\begin{remark}
It is also possible to use a  slightly different
groupoid: consider $\GR'=\tgt{G(M)}\times_{M\times M} \tgt{(U\times
U)}$. This is different from $\GR$ only at the level of $s=0, t>0$,
 which
does not matter as we will not use this part of the groupoid. For $s>0, t=0$, one has $A_M \times_{M\times M} (U\times
U)=A_M \times_M (U\times_M U)=A(M,U)$, and for $s>0, t>0$, $ G(M)\times_{M\times M} (U\times
U)=G(U)$. 
\end{remark}

Let us next denote by $e_{i,j}$ (for $i,j=0,1$) the various $K$-theory morphisms
induced by the restriction maps $\GR \to \GR_{s=i,t=j}$, where, for
instance, $\GR_{s=1,t=0} = A_M \times_M U \times_M U$ is the lower
right corner of the figure above.

\begin{lemma}
  The evaluation maps $e_{0, 0}$ and $e_{1, 0}$ are isomorphisms, and:
  \begin{itemize}
  \item  the analytic index of $U$ is   $\ind_a^U = e_{1, 1} \circ
    e_{0, 0}^{-1}$;
  \item the maps 
   \begin{multline*}
    \ind_a^{(1)} := e_{1, 0} \circ e_{0, 0}^{-1} : K^*(A^*_U)
    \to K^*(C^*(M, U)) \quad \text{and}\\
    \ind_a^{(2)} := e_{1, 1} \circ e_{1, 0}^{-1} :
    K^*(C^*(M, U)) \to K^*(C^*(U)).
\end{multline*}
are such that
\begin{equation}\label{eq.16}
    \ind_a^U = \ind_a^{(2)} \circ \ind_a^{(1)}.
\end{equation}
\end{itemize}
\end{lemma}

\begin{proof}
 We shall define several successive decompositions (figure \ref{fig:1}).
\begin{figure}[htbp]
  \centering
  \input{mn5.pstex_t}
  \caption{Successive restrictions of $\GR$}
  \label{fig:1}
\end{figure}
Let us notice that the above construction of $\GR$ is such that
$\GR' := \GR_1 \cup \GR_2$ is a closed subgroupoid of $\GR$, with
complement $\GR_3=G(U)\times (0,1]\times (0,1]$. This induces an exact
sequence
\begin{equation*}
0 \to C^*(\GR_3) \to C^*(\GR) \to C^*(\GR_1 \cup \GR_2) \to 0.
\end{equation*}
But $K_*(C^*(\GR_3))=0$, so that $K_*(C^*(\GR))\simeq K_*(C^*(\GR_1
\cup \GR_2))$.  Now $\GR_2$ is an open subgroupoid of $\GR_1 \cup
\GR_2$, with vanishing $K$-theory groups. A similar exact sequence
argument then shows that $K_*(C^*(\GR_1\cup \GR_2)) \simeq
K_*(C^*(\GR_1))\simeq K_*(A_U)$. Thus $e_{0, 0} : K_*(C^*(\GR)) \to
K_*(C_0(A_U^*))$ is an isomorphism. Since the restriction of $\GR$ to
the diagonal $s=t$ of $[0, 1]^2$ is the tangent groupoid of $U$, we
obtain that $\ind_a^U = e_{1, 1} \circ e_{0, 0}^{-1}$.

 We aim
at factorizing the index map $\ind_a^U$ in two maps defined
along the lower and right sides of the square. The first map is
\begin{equation*}
        \ind_a^{(1)} := e_{1, 0} \circ e_{0, 0}^{-1} : K^*(A^*_U) \to
        K_*(C^*(M, U)).
\end{equation*}
For the second map, we shall define other successive
  decompositions (figure \ref{fig:2}.)
\begin{figure}[htbp]
  \centering
  \input{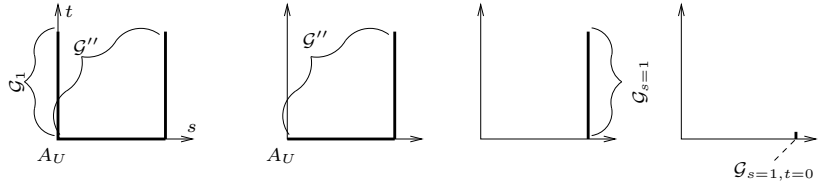}
  \caption{Other successive restrictions of $\GR$}
  \label{fig:2}
\end{figure}

For the second map, we consider the subgroupoid $\GR'':= \GR_{L}$
$L := \{s=1\}\cup \{t=0\}$. Using a short exact sequence argument as
above, we obtain that $\GR''\cup \GR_1$ has the same \kth\ as $\GR$
because its complement in $\GR$ is $G_U\times (0,1) \times (0,1]$.
Then, similarly, $\GR''\cup \GR_1$ turns out to have the same \kth\ as
$\GR''$ since the complement of $\GR''$ in $\GR''\cup \GR_1$ is $T_U
\times (0,1]$.
  Furthermore, the decomposition $\GR''=\GR_{s<1,t=0} \cup
\GR_{s=1}$ induces a $K$-isomorphism between $\GR''$ and
$\GR_{s=1}$. Indeed, still considering that $U=M\times \RR^n$, we
showed before that $\GR_{t=0}=A_M \times \tgt{G(\RR^n)}$. But the \kth\ of
$\tgt{G(\RR^n)}_{s<1}$ is isomorphic to that of the adiabatic groupoid,
which is zero since the analytic index for a euclidean space is an
isomorphism.

The last step is that the evaluation at $t=0$ in the groupoid
$\GR_{s=1}$ is an isomorphism, using the same arguments as before.
In conclusion, $e_{1,0}: \kc{\GR} \to K_*(C^*(M, U))$ is an isomorphism.

The map $\ind_a^{(2)} $ is thus well defined, and the equality
$$\ind_a^U = \ind_a^{(2)} \circ \ind_a^{(1)}$$
 is straightforward.
\end{proof}

\begin{remark}
  In the proof above we used $\GR''$ since we shall need later to
  consider the groupoid $\mathcal{H}=\GR_{\{s=1\}}$. But the proof could
  also have been handled considering the groupoid $\GR'$. 
\end{remark}

To end the proof of Proposition \ref{prop.diag1}, let us consider now the commutative diagram, in which the morphism
$\Theta_0$ is the isomorphism of the Equation \eqref{theta} (its
definition will be recalled below as part of a more general
construction)
\begin{equation}
\begin{CD}
    K^*(A_M^*) @>=>> K^*(A_M^*) @>{\ind_a^M}>> K_*(C^*(M))\\
    @V{k_!}VV  @V{\Theta_0}V{\simeq}V  @V{k_K}V{\simeq}V\\
    K^*(A_U^*) @>{\ind_a^{(1)}}>>
    K_*(C^*(M, U)) @>{\ind_a^{(2)}}>> K_*(C^*(U))
\end{CD}
\end{equation}
Equation \eqref{eq.16} shows that it is enough to prove that the
above diagram is commutative.

Let $X$ be the fiberwise one point compactification of $A_M^*$.
Then $K^*(A_M^*) \subset K^*(X)$, since $X \smallsetminus A_M^*$
is a retract of $X$. The commutativity of the left diagram then
follows from Corollary \ref{cor.family} (after we lift the bundle
$U$ to $X$).

We need to define the morphism $\Theta_0$. Let us consider the
groupoid $\maH$ defined as the restriction of $\GR$ to $\{1\}
\times [0, 1]$ used also to define the morphism $\ind_a^{(2)}$. It
has units $U \times [0, 1]$. Let $\Omega = r^{-1}(M \times [0,
1])$, as in the proof of Lemma \ref{lemma.prop.3}. As in the proof
of that lemma, $\Theta$ defines an imprimitivity module between
$\maH$ and $r^{-1}(M) \cap d^{-1}(M) = \tgt{\maG_M}$. This
imprimitivity module induces imprimitivity modules $\Theta_t$ for
$t \in [0, 1]$ (the parameter of the deformation). By the proof of
Proposition \ref{prop.3}, the isomorphism $k_K$ is defined by
$\Theta_1$. The isomorphism defined by $\Theta_0$ was also denoted
by $\Theta_0$. The commutativity of the right rectangle in the
above diagram then follows from the compatibility of the
isomorphisms defined by $\Theta$ with restriction morphisms.
\end{proof}

The commutativity of the second square in the Diagram \ref{diag2}
follows from the naturality of the tangent groupoid construction.
Here are the details.

\begin{proposition}\label{prop.diag2}\
Let $j : U \to X$ be the inclusion of the open subset $U$. Then
the diagram below commutes:
\begin{equation*}
\begin{CD}
    \kc{U} @>{j_K}>>\kc{X}\\
    @A{\ind_a^U}AA @AA{\ind_a^X}A\\
    K^*(A^*_U) @>{j_*}>\simeq> K^*(A^*_X).
\end{CD}
\end{equation*}
\end{proposition}

\begin{proof}
As $U$ is open in $X$, the \gr\ $\tgt{G(U)}$ identifies with the
restriction of  $\tgt{G(X)}$ to $U$. This induces a map $\tgt{\cc{G(U)}}
\to \cc{\tgt{G(X)}}$. So we get the following commutative diagram
\begin{equation*}
\begin{CD}
    \cc{U} @>>> \cc{X}\\
    @A{e_1}AA @AA{e_1}A\\
    \cc{\tgt{G(U)}} @>>> \cc{\tgt{G(X)}}\\
    @V{e_0}VV @VV{e_0}V\\
    \CO(A^*_U) @>>>  \CO(A^*_X),
\end{CD}
\end{equation*}
in which the vertical arrows are inclusions. This commutative
diagram, together with Lemma \ref{lemma.prop.3} and Remark
\ref{rem.prop.3}, give the analogue commutative diagram in
$K$-theory
\begin{equation*}
\begin{CD}
    \kc{G(U)} @>>> \kc{G(X)}\\
    @A{e_1}AA @AA{e_1}A\\
    K_*(C^*(\tgt{G(U)}))  @>>> K_*(C^*(\tgt{G(X)}))\\
    @V{e_0}VV @VV{e_0}V\\
    K^*(A^*(G(U))) @>>>  K^*(A^*(G(X)))\\
\end{CD}
\end{equation*}
The commutativity of the above diagram proves our result.
\end{proof}

We are ready now to prove one of our main results by putting
together what he have proved in the previous two propositions, as
explained in the beginning of this section.

\begin{theorem}\label{commutativity}
Let $M \overset{\iota}{\to} X$ be a closed embedding of manifolds with
corners. Then the  diagram
\begin{equation}\label{diag4}
\begin{CD}
    \kc{M)} @>{\iota_K}>> \kc{X}\\ @AA{\ind_a^M}A
    @A{\ind_a^X}AA\\ K^*(A^*_M) @>\iota_{!}>> K^*(A^*_X) \\
\end{CD}
\end{equation}
is commutative.
\end{theorem}

\begin{proof}
As $i(M)$ is a closed submanifold of $X$, there exists a tubular
neighborhood $U$ of $i(M)$ in $X$, along with a fibration $\pi: U
\to M$. Let $k : M \to U$ be the embedding of $M$ into $U$ as the
zero section and $j : U \to X$ be the embedding of $U$ as an open
subset of $M$. The result then follows from the commutativity of
the diagrams in Propositions \ref{prop.diag1} and \ref{prop.diag2}
and from $i_! = j_* \circ k_!$ and $i_K = j_K \circ k_K$. (The
diagram of Equation \eqref{diag2} explains this reasoning.)
\end{proof}

\section{An Atiyah--Singer type theorem}

Motivated by Theorem \ref{commutativity} and by the results of
Section \ref{sec.prop} (see Propositions \ref{prop.2} and
\ref{prop.3}) we introduce the following definition.

\begin{definition}\label{def.class}\
A classifying manifold $X_M$ of $M$ is a compact manifold with
corners $X_M$, together with a closed embedding $\iota: M \to X_M$
with the following properties:
\begin{enumerate}[(i)]
\item\ each open face of $X_M$ is diffeomorphic to a Euclidean space,
\item\ $F \to F \cap M$ induces a bijection between the open faces
of $X_M$ and ~$M$.
\end{enumerate}
\end{definition}

Note that if $M \subset X_M$ are as in the above definition, then
each face of $M$ is the transverse intersection of $M$ with a face
of $X_M$. As a consequence we obtain the following result, which
generalizes the main theorem of \cite{AS1}.

\begin{lemma}\label{lemma.t.i}\
Let $M$ be a \mec, and $\iota : M \hookrightarrow X_M$ be a
classifying space of $M$. Then the maps $\iota_K$ and $\ind_a^X$
of Theorem \ref{commutativity} are isomorphisms.
\end{lemma}

\begin{proof}
This was proved in Propositions \ref{prop.2} and \ref{prop.3}.
\end{proof}

Let $\iota : M \to X_M$ be a classifying space for $M$. The above
lemma then allows us to define (see the diagram \ref{diag4})
\begin{equation*}
    \ind_t^M :=\iota_K^{-1} \circ \ind_a^X \circ \iota_{!} :
    K^*(A^*_M) \to \kc{M}.
\end{equation*}

If $M$ is a smooth  compact manifold (so, in particular, $\pa M =
\emptyset$), then $C^*(M) = \maK$, the algebra of compact
operators on $L^2(M)$ and hence $K_0(C^*(M)) = \ZZ$. Any embedding
$\iota : M \hookrightarrow \RR^N$ will then be a classifying space for
$M$. Moreover, as explained in Remark \ref{rem.top}, for $X =
\RR^n$, the map $\iota_K^{-1} \circ \ind_a^X : K^*(TX) \to \ZZ$ is the
inverse of $j_! : K^0(pt) \to K^0(T\RR^N)$ and hence
$\ind_t^{\RR^N} = (j_!)^{-1} \iota_!$, which is the definition of the
topological index from \cite{AS1}. In view of this fact, we shall
also call the map $\ind_t^M$ {\em the topological index}
associated to $M$.

\begin{theorem}\label{thm.topological}\ The topological index
map $\ind_t^M$ depends only on $M$, that is, it is independent of
the classifying space $X_M$, and we have
\begin{equation*}
    \ind_t^M = \ind_a^M : K^*(A^*_M) \to \kc{M}.
\end{equation*}
\end{theorem}

\begin{proof}
This follows right away from Theorem \ref{commutativity}.
\end{proof}

If $M$ is a smooth compact manifold (without boundary), this
recovers the Atiyah-Singer index theorem on the equality of the
analytic and topological index \cite{AS1}.

\begin{remark} Let us also mention that $K_*(C^*(M)) \simeq
K^*(X_M)$ provides us with a way of determining $K_*(C^*(M))$,
which is a non-trivial problem.
\end{remark}

\section{Construction of the classifying space
$X_M$\label{construction}}

We now show that a classifying manifold $X_M$ of $M$ exists
(Definition \ref{def.class}). The choice of $X_M$ is not
canonical, in general.

Let $M$ be a compact \mec, and let $(H_i)_{1\leq i \leq r}$ be the
set of hyperfaces of $M$.  For each closed hyperface $H_i$, we
shall fix a defining function $\rho_i$ (these could be, for
example, the defining functions used in the definition of $G(M)$,
Definition \ref{def.GM}). Also, let us choose an embedding of
$\phi$ of $M$ into some $\R^N$. The map
\begin{equation*}
    \psi = (\phi, \rho_1, \ldots, \rho_r) : M \to
    \R^N \times [0, \infty)^r
\end{equation*}
is thus an embedding of manifolds with corners, which, however,
does not induce a bijection of the faces. To fix this problem, we
need to add extra coordinates which will disconnect the faces of
$\R^N\times [0, \infty)^r$ putting them in bijection with the
faces of $M$. If $J \subset \{1, \ldots, r\}$, define
\begin{equation*}
    F_J=\cap_{j \in J} H_j.
\end{equation*}
If nonempty, this is a disjoint union of closed faces of
codimension $|J|$, the number of elements of $J$. Assume $F_J$ is
not empty and let $f_J$ be a continuous function on $F_J$ with
values in $\{1, \ldots, n_J\}$ where $n_J$ is the number of
connected components of $F_J$. We chose $f_J$ to take different
values on different connected components of $F_J$. This function
can then be extended, thanks to Tietze's theorem, to a smooth
function still denoted by $f_J : M \to \R$. Let us denote by
$\maJ$ the set of nonempty subsets of $\{1, \ldots, r\}$ for which
$F_J$ is not empty, and $l=|\maJ|$ the number of its elements.
Then we obtain an embedding
\begin{equation*}
    \Psi := (\phi, \rho_1, \ldots, \rho_r, f_J)\, :\, M\,
    \to\, X_0 := \R^N\times [0, \infty)^r \times \R^l,
    \quad J \in \maJ.
\end{equation*}

We still need to disconnect the faces of $X_0$ whose inverse image
in $M$ is disconnected. To this end, for any $J \in \maJ$, let
$Y_J \subset \R^N\times [0, \infty)^r \times \R^l$ be the closed
subset defined by
\begin{equation*}
    Y_J := \{(z, x_1, \ldots, x_r, y_J),\
    y_J-\ha \in \ZZ\, \text{ and }\, x_j = 0\,
    \text{ for all }\, j \in J\}.
\end{equation*}
Then $\Psi(M) \cap Y_J = \emptyset$, by the construction of $f_J$.
Let $X_1 := X_0 \smallsetminus \cup_{J \in \maJ} Y_J$. Finally,
remove from $X_1$ all the faces that do not intersect $M$ and call
what is left $X_M$:
\begin{equation}\label{eq.def.XM}
    X_M := X_1 \smallsetminus \cup F,
    \qquad F \subset X_1 \, \text{ face such that }\,
    F \cap M = \emptyset.
\end{equation}

Naturally this creates many more faces than we have in $M$, so the
last step is to take $X$ to be the complementary in $X_1$ of the
open faces which do not intersect $\Psi(M)$.

\begin{proposition}
The manifold $X_M \supset M$ of Equation \eqref{eq.def.XM} is a
classifying space for $M$.
\end{proposition}

\begin{proof}\
We need to prove that
\begin{enumerate}
\item each open face of $X_M$ is diffeomorphic to a Euclidean space,
\item $F \to F \cap M$ defines a bijection between the set of open
faces of $M$ and the set of open faces of $X_M$,
\item $M$ is a closed, embedded submanifold of $X_M$.
\end{enumerate}

The {\em open} faces of $X_0 := \R^N\times [0, \infty)^r \times
\R^l$ are in one-to-one correspondence with the subsets of $\{1,
2, \ldots, r\}$. More precisely,
\begin{equation*}
    G_I = \{(z, x_1, \ldots, x_r, y_J)\in X_0,
    \ x_j = 0 \Leftrightarrow j \in I\},
\end{equation*}
$I \subset \{1, 2, \ldots, r\}$, are all the open faces of $X_0$.

Fix $I \subset \{1, 2, \ldots, r\}$. The open faces $F_1 \subset
X_1$ contained in $G_I$ are the connected components of
\begin{equation*}
    G_I \smallsetminus \cup_J Y_J =
    G_I \smallsetminus \cup_J (G_I \cap Y_J), \quad
    J \in \maJ.
\end{equation*}
Since $G_I \cap Y_J = \emptyset$ for any $J \in \maJ$ that is not
contained in $I$, it is enough to consider only $J \subset I$, $J
\in \maJ$. Fix a face $F_1 \subset X_1$. Then $y_J \in (m_J-1/2,
m_J+1/2)$ on $F_1$, for some $m_J \in \ZZ$. This shows that
\begin{equation*}
    F_1 = \{(z, \ldots , y_J) \in G_I,\
    m_J-1/2 < y_J < m_J+1/2,\ J \subset I, J \in \maJ \}.
\end{equation*}
In particular, $F_1 \simeq \RR^N \times (0, \infty)^a \times
(-1/2, 1/2)^b \times \RR^{l-b}$, $a=r-|I|$. This verifies
condition (a) above.

Let $F_1$ be the face fixed above. Then $F_1 \cap M$ is the set of
points $x \in M$ satisfying $\rho_i(x) = 0$, for all $i \in I$, and
$f_J(x) \in (m_J-1/2, m_J+1/2)$, for all $J \subset I$ such that $J
\in \maJ$, therefore
\begin{equation*}
    F_1 \cap M = \{ x \in M,\ \rho_i(x) = 0\, \text{ if } i \in I,\
    f_I(x) \in (m_I-1/2, m_I+1/2)\},
\end{equation*}
which is either empty or a connected component of $F_J$, by the
construction of $f_J$. In other words, for any open face $F_1 \subset
X_M$, the intersection $F_1 \cap M$ is either empty or an open face of
$M$. This verifies condition (b) above.

Finally, let us notice that each open face $F$ of $M$ is contained
in exactly one face $F_1$ of $X_M$ and each open face $F_1$ of
$X_1$ is contained in an open face $G_I$ of $X_0$. This shows that
$M$ intersects $F_1$ transversely, because $dx_j$, $j \in J$, are
linearly independent on $F_J := \cap_{j \in J} H_j$. Also, $M$ is
a closed submanifold of $X_M$ because $\phi : M \to \RR^N$ (the
first component of $\Psi$) is an embedding. This verifies
condition (c) above and thus completes the proof.
\end{proof}

\begin{remark}
In the case $M$ is a smooth manifold, our construction is
such that $X_M$ is a euclidean space.
\end{remark}

\bibliographystyle{amsplain}

\bibliography{mn}

\end{document}